\documentclass[11pt]{amsart}

\usepackage{color}
\usepackage{xspace}
\usepackage{enumerate}
\usepackage[utf8]{inputenc}

\usepackage{amssymb}
\usepackage{epsfig}
\usepackage{comment}
\usepackage{amsmath}
\usepackage{subfigure}
\usepackage[section]{placeins}
\usepackage{mathrsfs}
\usepackage{graphicx}
\usepackage[notrig]{physics}
\usepackage{units}
\usepackage{caption}
\usepackage{color}
\usepackage{pstricks}
\usepackage{esint}

\graphicspath
{
  {./}
    {./Figures/}
}

\begin{document}

\title{Homogenization and continuum limit of mechanical metamaterials}
\author{
M.~P.~Ariza${}^1$, S.~Conti${}^2$ and M. Ortiz${}^{2,3}$
}

\address
{
  ${}^1$Escuela T\'ecnica Superior de Ingenier\'ia, Universidad de Sevilla,
  Camino de los descubrimientos, s.n., 41092 Sevilla, Spain.
  \\
   ${}^2$Institut f\"ur Angewandte Mathematik and Hausdorff Center for Mathematics, Universit\"at Bonn, Endenicher Allee 60, 53115 Bonn, Germany.
  \\
   ${}^3$Division of Engineering and Applied Science, California Institute of Technology, 1200 E.~California Blvd., Pasadena, CA 91125, USA.
}

\email{$\dots$}

\begin{abstract}
When used in bulk applications, mechanical metamaterials set forth a multiscale problem with many orders of magnitude in scale separation between the micro and macro scales. However, mechanical metamaterials fall outside conventional homogenization theory on account of the flexural, or bending, response of their members, including torsion. We show that homogenization theory, based on calculus of variations and notions of Gamma-convergence, can be extended to account for bending. The resulting homogenized metamaterials exhibit intrinsic generalized elasticity in the continuum limit. We illustrate these properties in specific examples including two-dimensional honeycomb and three-dimensional octet-truss metamaterials.
\end{abstract}

\maketitle

\begin{center}
Dedicated to Alan Needleman, {\sl ab imo pectore}, on the occasion of his 80th birthday.
\end{center}

\section{Introduction}

{\sl Mechanical metamaterials} \cite{Montemayor2015, Greer:2019, LU2022}, or {\sl architected materials}, are characterized by a reticular structure that affords great flexibility for tailoring and shaping their mechanical properties (cf., e.~g., \cite{Fleck2010}, for a review of stiffness, strength and damage tolerance; \cite{BENEDETTI2021}, for a review of fatigue performance; and \cite{WU2021}, for a review of dynamic properties). Advances in fabrication techniques \cite{Greer:2019, Xia2022, Jin2024, HAMZEHEI2024} have brought within the realm of practicality new metamaterial designs that populate previously inaccessible areas of material-property space \cite{Ashby2006}, including lightweight structures with high strength and elastic modulus \cite{lee2012a, jang2013a, bauer2014a, meza2014a, rys2014a, zheng2014a, meza2015a, gu2015a, bauer2015a, Montemayor2015, bauer2016a}. Common designs include honeycomb structures \cite{bauer2014a, bauer2016a} and the octet-truss \cite{deshpande2001b, jang2013a, bauer2014a,  zheng2014a, meza2015a, gu2015a, bauer2015a}, among others \cite{lee2012a, bauer2014a,  rys2014a, ros2015a, bauer2016a, Rosario2017}. Metamaterials often exhibit scaling properties and size effects \cite{deshpande2001a, deshpande2001b} that exploit the ‘‘smaller is stronger’’ principle to beneficial effect, especially as regards fracture resistance \cite{CHEN1998, Shaikeea2022, Maurizi2022, Luan2022} and energy absorption \cite{ZHANG2023}.

\begin{figure}[h]
\begin{center}
    \includegraphics[width=8cm]{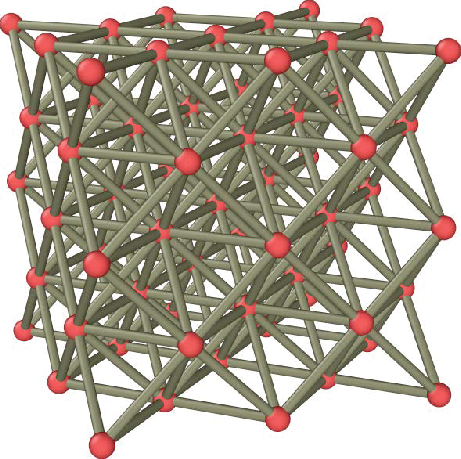}
\end{center}
\caption{Octet-truss metamaterial \cite{deshpande2001a, deshpande2001b}.}
\label{bCb5kk}
\end{figure}

Mechanical nanostructured metamaterials, when used in bulk applications, define a quintessential multiscale problem, with many orders of magnitude in scale separation between the micro and macro domains. Under such conditions, models---computational or otherwise---predicated on the explicit tracking of every individual member of the metamaterial are impractical and uncalled for, since locally, away from defects and stress concentrations, large numbers of members deform collectively under applied loading and jointly exhibit well-defined effective behavior that can be characterized by means of a continuum material law. Coarse-graining approaches \cite{PHLIPOT2019, kochmann:2019} based on the quasicontinuum method \cite{Tadmor:1996a, Tadmor:1996b, Espanol:2013}, data-driven approaches \cite{Korzeniowski:2022, Weinberg:2023, Jiao2023} and machine-learning approaches \cite{Karapiperis2023} have been proposed in order to bridge length scales while providing full-resolution of the lattice where necessary.

In this work, we seek to characterize analytically the effective behavior of metamaterials at the macroscale under conditions of strict separation between the structural and lattice scales. The identification of such effective material behavior is the aim of homogenization theory \cite{Cioranescu:1999} and discrete-to-continuum techniques \cite{Cicalese:2009, BraidesGelli2006}. However, mechanical metamaterials fall outside conventional homogenization theory on account of the flexural, or bending, response of their members.

We show that homogenization theory can be extended to account for bending and that the resulting homogenized metamaterials exhibit generalized elasticity of the {\sl micropolar} type \cite{Eringen:1966}. Such theories are often used to introduce an {\sl ad hoc} micromechanical length into the material law (cf., e.~g., \cite{Askes:2011} for a review and discussion). By contrast, in the present context the micropolar terms are physical and are introduced by bending. In addition, the form of the limiting continuum energy is identified uniquely and unequivocally by the homogenization analysis. We additionally show that, for linear-elastic matematerials, the effective moduli of the limiting continuum energy can be computed explicitly using the discrete Fourier transform. We illustrate these properties in specific examples including two-dimensional honeycomb and three-dimensional octet-truss metamaterials.

\section{Mechanical metamaterials}
\label{7DX37E}

We consider throughout {\sl open-cell mechanical metamaterials} comprised of bars, possessing both axial and bending stiffness, connected at rigid joints capable of transmitting forces and moments.

\subsection{Notational conventions}

We denote by $\mathbb{Z}$ the set of integers; by $\mathbb{R}$ the set of real numbers; and by $\mathbb{C}$ the set of complex numbers.

The {\sl characteristic function} of a subset $A\subset X$ is defined as
\begin{equation} \label{qe7B0S}
    \chi_A(x)
    =
    \left\{
        \begin{array}{ll}
            1, & \text{if} \;\; x \in A , \\
            0, & \text{if} \;\; x \in X\backslash A .
        \end{array}
    \right.
\end{equation}
For a Lebesgue-measurable set $A \subset \mathbb{R}^n$, we denote by $| A |$ its Lebesgue measure, i.~e., its length if $n=1$; its area if $n=2$; and its volume if $n=3$.

We denote by $x\cdot y \in \mathbb{R}$ the {\sl dot product} of two vectors $x$, $y\in\mathbb{R}^n$; by $y \otimes x \in \mathbb{R}^{n\times m}$ the {\sl tensor product} of two vectors $x \in \mathbb{R}^m$ and $y \in \mathbb{R}^n$, i.~e., $(y \otimes x)_{ji} = y_j \, x_i$; and by $y=A x \in \mathbb{R}^n$ the product between a matrix $A \in \mathbb{R}^{n\times m}$ and a vector $x \in \mathbb{R}^m$. We use throughout the abbreviation
\begin{equation}
    A x \cdot y := (A x) \cdot y = \sum_{j=1}^n \sum_{i=1}^m A_{ji} y_j x_i ,
\end{equation}
corresponding to assigning precedence to matrix multiplication over dot product. For a matrix $A \in \mathbb{R}^{n\times n}$ we denote by $A^T$ its transpose and by
\begin{equation}
    \operatorname{\rm sym} A = \frac{1}{2} (A+A^T) ,
    \quad
    \operatorname{\rm skw} A = \frac{1}{2} (A-A^T) ,
\end{equation}
its symmetric and skewsymmetric parts.

Given two vectors $x \in \mathbb{R}^m$ and $y \in \mathbb{R}^n$, we denote by $z = (x;y) \in \mathbb{R}^{m}\times\mathbb{R}^{n}$ the {\sl cartesian product} of $x$ and $y$, i.~e., $z_i = x_i$, for $i=1,\dots,m$ and $z_{j+m} = y_j$, for $j=1,\dots,n$. Given two points $x$, $y \in \mathbb{R}^n$, we denote by $[x,\,y]$ the {\sl segment} of straight line joining the two points in $\mathbb{R}^n$.

Given a function $f:\mathbb{R}^n \to \mathbb{R}^m$, we denote by $\partial_i f(x)$ its partial derivative with respect to $x_i$. For functions of one variable, we also write $f'(x) = \partial f(x)$. We denote by $Df : \mathbb{R}^n \to \mathbb{R}^{m\times n}$ the matrix of partial derivatives of $f$, i.~e., $(Df)_{ij} = \partial_j f_i$. Given a Lebesgue-measurable subset $A \subset \mathbb{R}^n$
and an integrable function $f:A\to\mathbb R$,
we denote by
\begin{equation}
    I = \int_A f(x) \, dx
\end{equation}
the Lebesgue integral of $f(x)$ over $A$.

\subsection{Geometry of metamaterials}
\label{seckinematics}

Mechanical metamaterials are locally periodic frame structures consisting of {\sl bars} connected at a set of {\sl joints} (cf., e.~g., \cite{Fleck2010, Montemayor2015}. {\sl Infinite} metamaterials are comprised of $M$ types, or {\sl classes}, of oriented bars, all the bars in a class being identical, including orientation, modulo translations. For instance, the bar classes of a honeycomb metamaterial are $\{\diagdown,\, \text{---},\, \diagup \}$, or $M=3$, Fig.~\ref{FIDpWS}. The {\sl local environment} of a joint is the set of bars incident on the joint. Metamaterials contain $N$ different types, or classes, of local environments, all local environments in a class being identical modulo translations. Correspondingly, the joints of the metamaterial are classified into $N$ types, according to their local environment. For instance, the joint classes of a honeycomb metamaterial are $\{\rotatebox[origin=c]{90}{${\textsf{Y}}$},\, \rotatebox[origin=c]{-90}{${\textsf{Y}}$}\}$, or $N=2$, Fig.~\ref{FIDpWS}.

The joints and bars of an infinite metamaterial are embedded periodically in space so that joints and bars of the same class span shifted {\sl Bravais lattices}. Thus, the positions of the joints in each class span a point set of the form
\begin{equation} \label{M6JjBD}
  x(l,\alpha) = b_\alpha + \sum_{i=1}^n l^i a_i ,
  \quad l \in \mathbb{Z}^n , \quad \alpha=1,\dots,N ,
\end{equation}
where $n$ is the dimension of space, e.~g., $n=2$ for two-dimensional lattices and $n=3$ for three-dimensional lattices; $(a_i)_{i=1}^n$ is a {\sl basis} of $\mathbb{R}^n$; $l := (l^i)_{i=1}^n$ is a corresponding array of integer {\sl lattice coordinates}; and $(b_\alpha)_{\alpha=1}^N$ are translation vectors, or {\sl shifts}. We additionally denote by $V$ the volume of the unit cell of the spanning Bravais lattice of the metamaterial. The joints of a metamaterial can therefore be {\sl indexed} by a pair $(l,\alpha)$, where $\alpha \in \{1,\dots,N\}$ designates the joint class and $l \in \mathbb{Z}^n$ its lattice coordinate array.

A full description of an infinite metamaterial additionally requires a scheme for indexing its bars by type and position, and for tabulating their connections to joints. By periodicity and translation invariance, the bars of a given bar class span the same Bravais lattice as the joint classes, modulo translations. Thus, bars can be indexed by pairs $(m,\beta)$, where $\beta \in\{ 1,\dots,M \}$ designates the bar class and $m \in \mathbb{Z}^n$ is a lattice coordinate array. In addition, an {\sl orientation} of the bar classes identifies beginning and end joints in every bar, which we designate by the symbols $\pm$, respectively (cf.~Section~\ref{vX3HxB} for examples).

The adjacency relations between bars and joints, or {\sl connectivity}, can be defined by designating, for every bar $(m,\beta)$, the corresponding beginning and end joints, $(m^-,\beta^-)$ and $(m^+,\beta^+)$, respectively. We note that, by periodicity,
\begin{equation} \label{g1CMLJ}
    ((m+t)^\pm,\beta^\pm) = (m^\pm + t,\beta^\pm),
    \quad
    \text{for all} \;\; t \in \mathbb{Z}^n ,
\end{equation}
i.~e., the adjacency relations are {\sl translation-invariant}. The corresponding joint coordinates follow from (\ref{M6JjBD}) as
\begin{equation} \label{46x3dB}
    x^\pm(m,\beta)
    :=
    x(m^\pm,\beta^\pm)
    =
    b_{\beta^\pm} + \sum_{i=1}^n (m^\pm)^i a_i ,
\end{equation}
and the {\sl spanning vector} of the joints follows as
\begin{equation}
    dx_\beta
    :=
    x^+(m,\beta) - x^-(m,\beta) ,
\end{equation}
which is independent of $m$ by the translation invariance property (\ref{g1CMLJ}). We also denote by
\begin{equation} \label{AIt7cu}
    h(m,\beta)
    =
    [x^-(m,\beta),\, x^+(m,\beta)]
\end{equation}
the spanning segments of the bars, regarded as sets.

We refer to subsets of the infinite metamaterials just described as {\sl metastructures}. The extent of a metastructure, i.~e., its set of bars, may be characterized by means of an index set
\begin{equation} \label{0Jdwod}
\begin{split}
    &
    J
    =
    \{ (m,\beta) \in \mathbb{Z}^n\times \{1,\dots, M\} \, : \,
    \\ & \qquad\qquad
    \text{bar} \; (m,\beta) \; \text{in metastructure} \} .
\end{split}
\end{equation}
For an infinite metamaterial, $J = \mathbb{Z}^n\times \{1,\dots, M\}$.

\subsection{Elastic energy of metastructures} \label{7PHmId}

Within the framework of linearized kinematics, the joints of a metastructure are endowed with deflection and rotation-angle degrees of freedom, $v(l,\alpha) \in \mathbb{R}^n$ and $\theta(l,\alpha) \in \mathbb{R}^{n(n-1)/2}$, respectively. We further denote by
\begin{equation} \label{FFPJsh}
  u(l,\alpha) := (v(l,\alpha); \theta(l,\alpha))
\end{equation}
the set of degrees of freedom of joint $(l,\alpha)$. At the local level, we denote by
\begin{equation} \label{0MjvOi}
  {U}(m,\beta) := (u^-(m,\beta); u^+(m,\beta))
\end{equation}
the set of local degrees of freedom of bar $(m,\beta)$, where $u^\pm(m,\beta)$ are the degree-of-freedom arrays of the joints of the bar.

The total energy of a linear-elastic metastructure of extent $J$ has the form
\begin{equation} \label{1ElWKO}
    E(u)
    =
    \sum_{(m,\beta)\in J}
    \frac{1}{2} \,
    \mathbb{S}_\beta \, {U}(m,\beta) \cdot {U}(m,\beta),
\end{equation}
where $\mathbb{S}_\beta$ is the elastic stiffness of a bar of class $\beta$. Throughout this work, we consider straight bars of constant cross section. Other cases, such as bars in the form of curved beams, tapered beams, and others, can be treated likewise.

\subsection{Straight two-dimensional beam of constant cross section} \label{yPlScU}

\begin{figure}[h]
\begin{center}
    \includegraphics[width=0.4\textwidth]{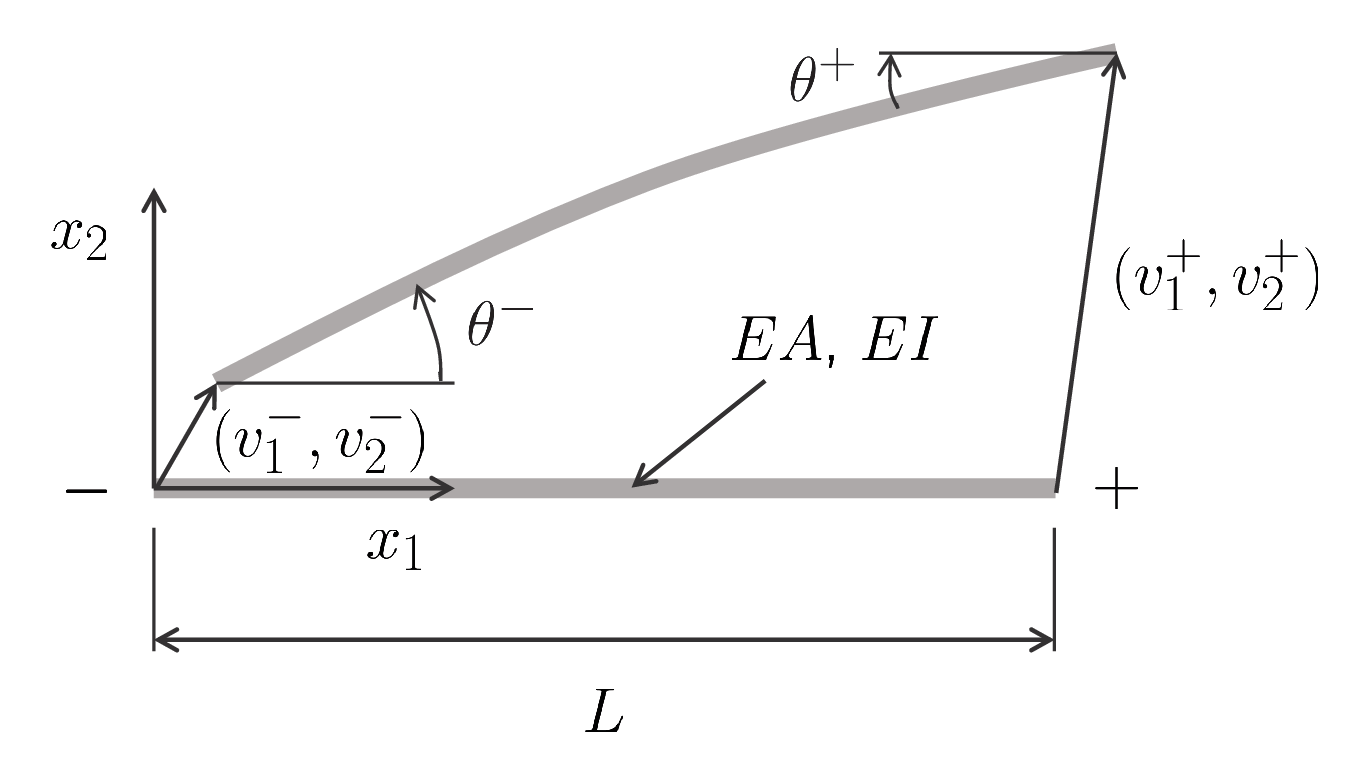}
\end{center}
\caption{Straight two-dimensional beam of constant cross section. Schematic of reference configuration and local deformation degrees of freedom.}
\label{DiuMhe}
\end{figure}

Two-di\-men\-sion\-al metamaterials are characterized by planar geometries with joint coordinates $x(l,\alpha) \in \mathbb{R}^2$, $n=2$. In addition, we assume that the metamaterial remains planar after deformation, so that the joint degrees of freedom $u(l,\alpha) := (v(l,\alpha);\, \theta(l,\alpha))$ of the joints comprise in-plane deflections $v(l,\alpha) \in \mathbb{R}^2$ and a single rotation angle $\theta(l,\alpha) \in \mathbb{R}$.

In order to represent the energy of the metamaterial, we consider a straight beam of length $L$ and constant cross section deforming in the $(x_1,x_2)$-plane, in a reference configuration spanning the domain $[0,L]$ on the $x_1$ axis, Fig.~\ref{DiuMhe}. We denote by $u^- := (v^-; \theta^-)$ and $u^+ := (v^+; \theta^+)$ the degrees of freedom of the beginning and end joint of the beam respectively. The corresponding linearized rotation matrices are
\begin{equation} \label{RZBy8m}
    w^\pm
    =
    \left(
        \begin{array}{cc}
            0 & - \theta^\pm \\ \theta^\pm & 0
        \end{array}
    \right)
    =:
    *\theta^\pm ,
    \quad
    \theta^\pm = *w^\pm ,
\end{equation}
where $*$ denotes the Hodge-* operator. We additionally denote by ${U} = (u^-; u^+)$ the array of local degrees of freedom.

In this representation, an application of elementary beam theory \cite{Timoshenko:1965} gives the energy of the beam as
\begin{equation} \label{7HCHbN}
\begin{split}
    E_{\rm ref}({U})
    & =
    \frac{{EA} }{2 L} \, ({v_1^+}-{v_1^-})^2
    +
    \frac{{EI} }{2 L} \, ({\theta^+}-{\theta^-})^2
    \\ & +
    \frac{6 {EI}}{L}
    \Big(
        \frac{{v_2^+}-{v_2^-}}{L}
        -
        \frac{{\theta^-} +{\theta^+}}{2}
    \Big)^2,
\end{split}
\end{equation}
where $E$ is the Young's modulus of the beam, $A$ the cross-sectional area, $I$ the moment of inertia of the cross section, with principal axis of inertia perpendicular to the plane. The first term in (\ref{7HCHbN}) is the energy due to the axial deformation of the beam, whereas the second and third terms account for the energy due to bending.

The quadratic energy (\ref{7HCHbN}) can also be expressed in matrix form in terms of a $6\times 6$ stiffness matrix $\mathbb{S}_{\rm ref}$, whence the stiffness matrix $\mathbb{S}_\beta$ of the bars of class $\beta$ in (\ref{1ElWKO}) follows as
\begin{equation} \label{QRbyOM}
    \mathbb{S}_\beta = T_\beta^T \, \mathbb{S}_{\rm ref} \, T_\beta ,
\end{equation}
in terms of some suitable, class-dependent, orthonormal change-of-axes matrix $T_\beta$.

By direct inspection of (\ref{7HCHbN}), we verify that the global energy (\ref{1ElWKO}) is {\sl material-frame indifferent}, i.~e., invariant under  transformations of the form
\begin{subequations} \label{UdCT8i}
\begin{align}
    &
    v(l,\alpha) \mapsto v(l,\alpha) + w \, x(l,\alpha) + c ,
    \\ &
    \theta(l,\alpha) \mapsto \theta(l,\alpha) + *w,
\end{align}
\end{subequations}
where $w \in so(2)$ represents a rigid linearized rotation and $c \in \mathbb{R}^2$ a rigid translation in the plane.

From these properties, the global energy (\ref{1ElWKO}) follows directly in global coordinates as
\begin{equation} \label{JqkpoV}
\begin{split}
    E(u)
    & =
    \sum_{(m,\beta) \in J}
    \Big\{
        \frac{{EA}_\beta }{2 L_\beta}
        \Big(
            dv(m,\beta) \cdot d_{1,\beta}
        \Big)^2
        \\ & +
        \frac{{EI}_\beta }{2 L_\beta} \,
        d\theta^2(m,\beta)
        \\ & +
        \frac{6 {EI}_\beta}{L_\beta}
        \Big(
            \frac{dv(m,\beta) \cdot d_{2,\beta}}{L_\beta}
            -
            \bar{\theta}(m,\beta)
        \Big)^2
    \Big\} ,
\end{split}
\end{equation}
where we write
\begin{subequations} \label{moEHtz}
\begin{align}
    &   \label{9JzElM}
    dv(m,\beta) = v(m^+,\beta^+)-v(m^-,\beta^-) ,
    \\ &    \label{cij5hA}
    d\theta(m,\beta) = \theta(m^+,\beta^+)-\theta(m^-,\beta^-) ,
    \\ &    \label{Fg7qSy}
    \bar{\theta}(m,\beta)
    =
    \frac{\theta(m^+,\beta^+)+\theta(m^-,\beta^-)}{2} ,
\end{align}
\end{subequations}
and
\begin{equation} \label{ARJlKT}
    d_{1,\beta} := \frac{1}{L_\beta} \, dx_\beta ,
    \quad
    d_{2,\beta} := d_{1,\beta}^\perp ,
\end{equation}
are orthonormal {\sl director}s for each bar class, with $d_{1,\beta}$ a unit vector aligned with the oriented axis and $d_{2,\beta}$ a unit vector perpendicular to $d_{1,\beta}$ with orientation chosen such that $(d_{1,\beta}, d_{2,\beta})$ defines a right-handed orthonormal basis. The variables
\begin{subequations}
\begin{align}
    &
    \epsilon_\beta = \frac{1}{L_\beta} dv(m,\beta)\cdot d_{1,\beta} ,
    \\ &
    \chi_\beta = \frac{1}{L_\beta} d\theta(m,\beta) ,
\end{align}
\end{subequations}
in (\ref{JqkpoV}) may be interpreted as local {\sl axial} and {\sl bending strains}, the first and second terms in (\ref{JqkpoV}) then supplying the corresponding axial and bending energies, respectively.

\subsection{Straight three-dimensional beam of constant cross section} \label{Akno0S}

\begin{figure}[h]
\begin{center}
    \subfigure[]{\includegraphics[width=0.4\textwidth]{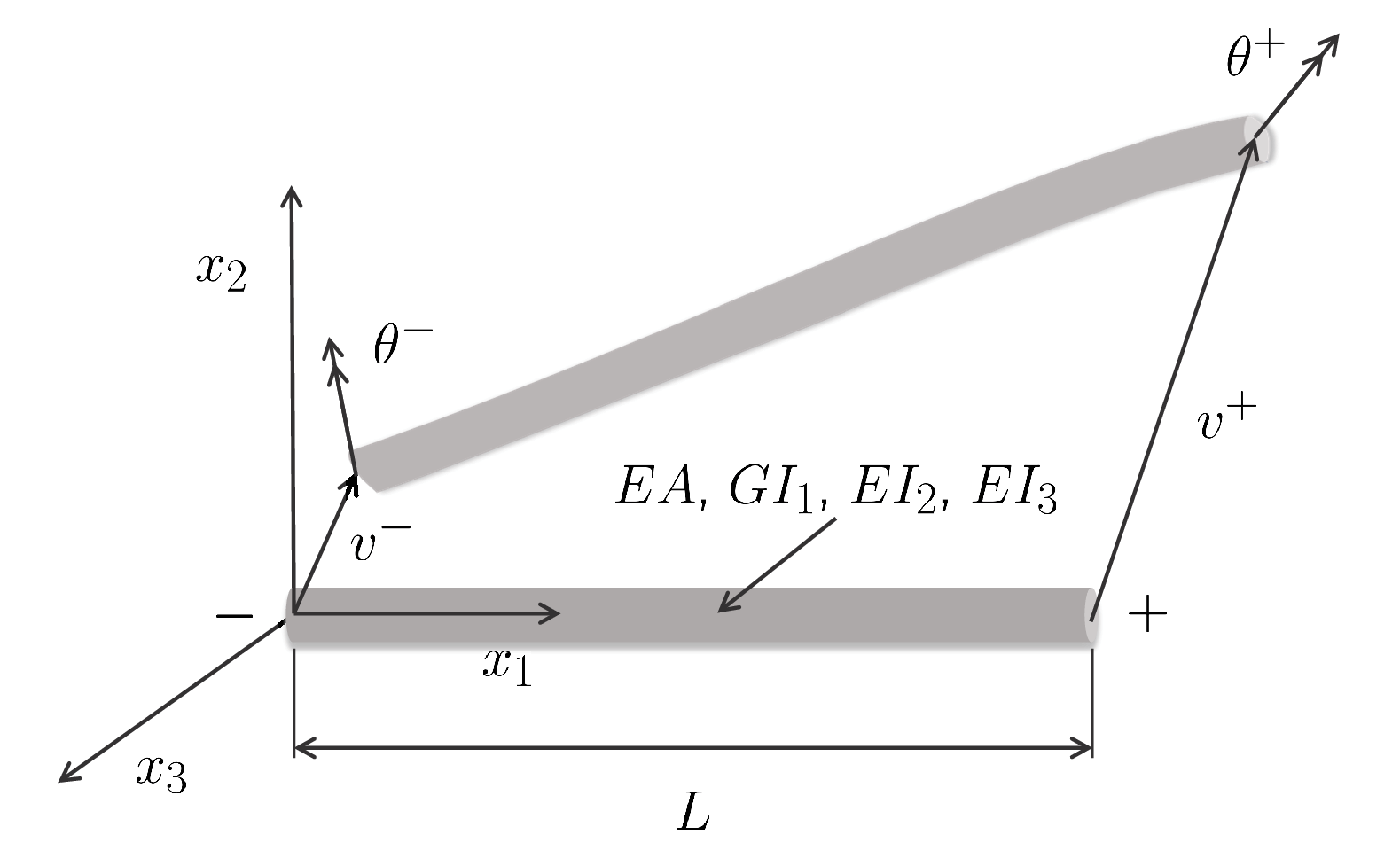}}
   \subfigure[]{\includegraphics[width=0.125\textwidth]{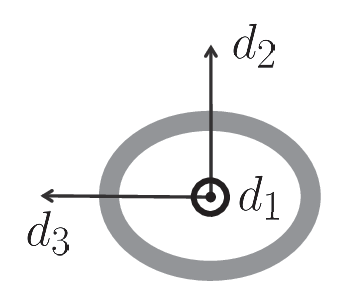}}
\end{center}
\caption{Straight three-dimensional beam of constant cross section. a) Schematic of reference configuration and local deformation degrees of freedom. b) Cross section and directors.}
\label{gWRlSr}
\end{figure}

Three-dimensional metamaterials are characterized by spatial geometries with joint coordinates $x(l,\alpha) \in \mathbb{R}^3$, $n=3$, and joint degrees of freedom $u(l,\alpha) := (v(l,\alpha);\, \theta(l,\alpha))$ comprising deflections $v(l,\alpha) \in \mathbb{R}^3$ and rotation angles $\theta(l,\alpha) \in \mathbb{R}^3$.
With this parametrization, the corresponding linearized rotation matrices are
\begin{subequations} \label{dkANTL}
\begin{align}
    &
    w^\pm
    =
    \left(
        \begin{array}{ccc}
            0 & -\theta_3^\pm & \;\;\theta_2^\pm \\
            \;\;\theta_3^\pm & 0 & -\theta_1^\pm \\
            -\theta_2^\pm & \;\;\theta_1^\pm & 0
        \end{array}
    \right)
    :=
    *\theta^\pm,
    \\ &
    \theta^\pm = *w^\pm,
\end{align}
\end{subequations}
where we use the Hodge-* operator. With ${U} = (u^-; u^+)$ as the array of local degrees of freedom, elementary beam theory \cite{Timoshenko:1965} now gives the energy of the beam as
\begin{equation} \label{tuxI9g}
\begin{split}
    E_{\rm ref}({U})
    & =
    \frac{{EA} }{2 L} \, ({v_1^+}-{v_1^-})^2
    + 
    \frac{{GI_1} }{2 L} \, ({\theta_1^+}-{\theta_1^-})^2
    \\ & +
    \frac{EI_2}{2 L} \, ({\theta_2^+}-{\theta_2^-})^2
    + 
    \frac{{EI_3} }{2 L} \, ({\theta_3^+}-{\theta_3^-})^2
    \\ & +
    \frac{6 EI_2}{L}
    \Big(
        \frac{v_3^+ - v_3^-}{L}
        +
        \frac{\theta_2^-+\theta_2^+}{2}
    \Big)^2
    \\ & +
    \frac{6 EI_3}{L}
    \Big(
        \frac{v_2^+ - v_2^-}{L}
        -
        \frac{\theta_3^-+\theta_3^+}{2}
    \Big)^2,
\end{split}
\end{equation}
where $A$ its cross-sectional area, $I_1$ its torsional moment of inertia, $I_2$ and $I_3$ its principal moments of inertia about axes $x_2$ and $x_3$, respectively, $E$ is the Young's modulus and $G$ is the shear modulus.

As in the two-dimensional case, the quadratic energy (\ref{tuxI9g}) can be expressed in matrix form in terms of a $12\times 12$ stiffness matrix $\mathbb{S}_{\rm ref}$, and the stiffness matrices $\mathbb{S}_\beta$ of the bars of class $\beta$ in cf.~(\ref{1ElWKO}) then follow as in (\ref{QRbyOM}) in terms of orthonormal transformation matrices $T_\beta$. The resulting global energy (\ref{1ElWKO}) is again material-frame-indifferent, i.~e., invariant under  transformations of the form (\ref{UdCT8i}),
with $w \in so(3)$ now representing a rigid linearized rotation and $c \in \mathbb{R}^3$ a rigid translation in three-dimensional space.

Alternatively, from these properties the global energy (\ref{1ElWKO}) follows directly in global coordinates as
\begin{equation} \label{w1G0fN}
\begin{split}
    &
    E(u)
    =
    \sum_{(m,\beta) \in J}
    \Big\{
        \frac{{EA}_\beta }{2 L_\beta}
        \Big(
            dv(m,\beta) \cdot d_{1,\beta}
        \Big)^2
        + \\ & \quad
        \frac{{GI}_{1,\beta} }{2 L_\beta} \,
        \Big(
            d\theta(m,\beta) \cdot d_{1,\beta}
        \Big)^2
        + \\ & \quad
        \frac{{EI}_{2,\beta} }{2 L_\beta} \,
        \Big(
            d\theta(m,\beta) \cdot d_{2,\beta}
        \Big)^2
        + \\ & \quad
        \frac{{EI}_{3,\beta} }{2 L_\beta} \,
        \Big(
            d\theta(m,\beta) \cdot d_{3,\beta}
        \Big)^2
        + \\ & \quad
        \frac{6 {EI}_{2,\beta}}{L_\beta}
        \Big(
            \frac{dv(m,\beta) \cdot d_{3,\beta}}{L_\beta}
            +
            \bar{\theta}(m,\beta) \cdot d_{2,\beta}
        \Big)^2
        + \\ & \quad
        \frac{6 {EI}_{3,\beta}}{L_\beta}
        \Big(
            \frac{dv(m,\beta) \cdot d_{2,\beta}}{L_\beta}
            -
            \bar{\theta}(m,\beta) \cdot d_{3,\beta}
        \Big)^2
    \Big\}
\end{split}
\end{equation}
where $dv(m,\beta)$, $d\theta(m,\beta)$ and $\bar{\theta}(m,\beta) $ are as in (\ref{moEHtz}) and $(d_{1,\beta}, d_{2,\beta}, d_{3,\beta})$ is an orthonormal {\sl director triad} in which $d_{1,\beta}$ is the direction of the bar axis, eq.~(\ref{ARJlKT}), and $d_{2,\beta}$ and $d_{3,\beta}$ are the principal directions of inertia of the cross section.

As in the two-dimensional case, the variables
\begin{subequations}
\begin{align}
    &
    \epsilon_\beta = \frac{1}{L_\beta} dv(m,\beta)\cdot d_{1,\beta} ,
    \\ &
    \chi_{i,\beta} = \frac{1}{L_\beta} d\theta_i(m,\beta) ,
\end{align}
\end{subequations}
may be interpreted as local {\sl axial} and {\sl bending strains}, the axial and bending energies in (\ref{w1G0fN}) then being quadratic in the corresponding strains.

\section{The continuum limit of metamaterials}

Next, we concern ourselves with the limit in which the metamaterial cell-size is much smaller than the domain size, or {\sl continuum limit}. To this end, we appeal to tools from calculus of variations, especially methods of $\Gamma$-convergence \cite{dalmaso:1993}, as they bear on discrete-to-continuum problems \cite{Alicandro:2004,BraidesGelli2006}. The fundamental property of variational methods is that they ensure convergence of the solutions of the finite cell-size problem to the solutions of the limiting continuum model. However, we note that the local energy of the bars in (\ref{JqkpoV}) and (\ref{w1G0fN}) includes terms that penalize the difference between the deflection-compatible rotations and the mean rotations  of the bars, e.~g., $dv(m,\beta) \cdot d_{2,\beta} / L_\beta$ and $\bar{\theta}(m,\beta) $ in (\ref{JqkpoV}) and similar terms in (\ref{w1G0fN}), respectively. This local energy structure renders the attendant homogenization analysis non-standard and outside the scope of conventional discrete-to-continuum results (cf.~\cite{Alicandro:2004,BraidesGelli2006}.

In order to characterize the continuum limit, we consider throughout metastructures whose extent is defined by a domain $\Omega\subset \mathbb{R}^n$. Specifically, the metastructure {\sl maximally contained} in $\Omega$ is defined by the index set
\begin{equation} \label{ue5kdD}
    J_\Omega
    =
    \{
        (m,\beta) \in \mathbb{Z}^n\times \{1,\dots, M\}\, : \,
        h(m,\beta)\subset \Omega
    \}
\end{equation}
with $x^\pm(m,\beta)$ as in (\ref{46x3dB}) and $h(m,\beta)$ as in (\ref{AIt7cu}). Thus, (\ref{ue5kdD}) records all the bars that are fully contained in $\Omega$. It also supplies a criterion for {\sl terminating} the metastructure at the boundary. We additionally denote by $I_\Omega$ the index set of\emph{} the joints in the metastructure, i.~e.,
\begin{equation}
\begin{split}
    &
    I_\Omega
    =
    \{
        (l,\alpha) \in \mathbb{Z}^n\times \{1,\dots, N\} \, : \,
        \operatorname{\exists} (m,\beta) \in J_\Omega
        \\ & \qquad
        \text{s.~t.} \;\;
        (l,\alpha) = (m^-,\beta^-)
        \;\; \text{or} \;\;
        (l,\alpha) = (m^+,\beta^+)
    \} .
\end{split}
\end{equation}

\begin{figure*}[ht!]
\begin{center}
    \includegraphics[width=13cm]{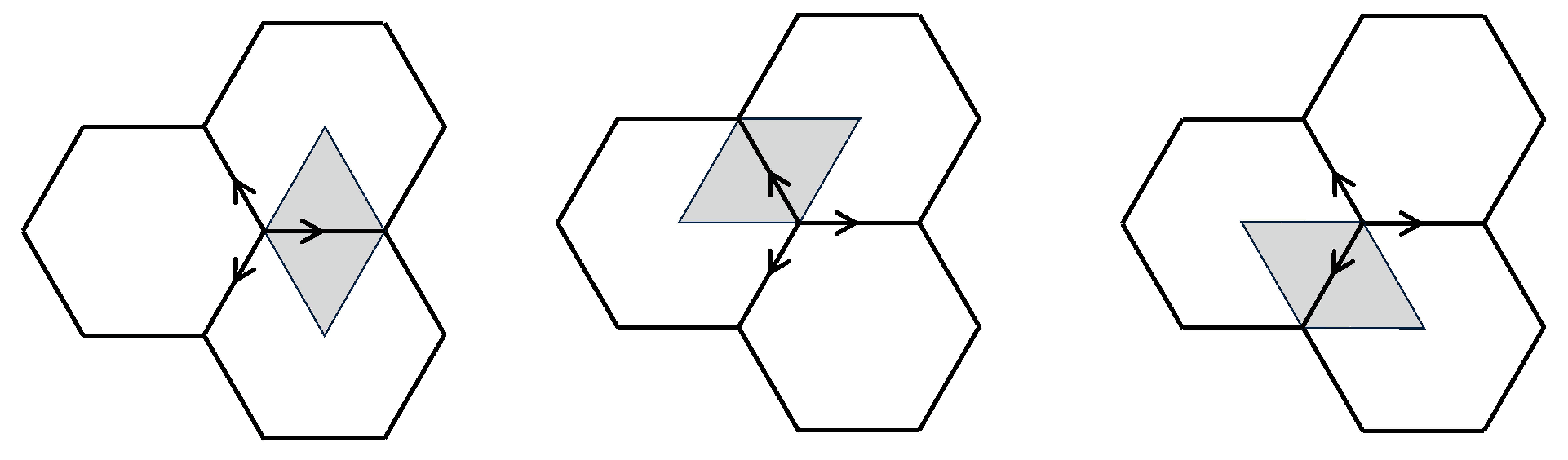}
\end{center}
\caption{Bar-wise Voronoi cells of a honeycomb metamaterial relative to its hexagonal cells.}
\label{bCb5vv}
\end{figure*}

\begin{figure*}
\begin{center}
    \subfigure[]{\includegraphics[width=0.4\textwidth]{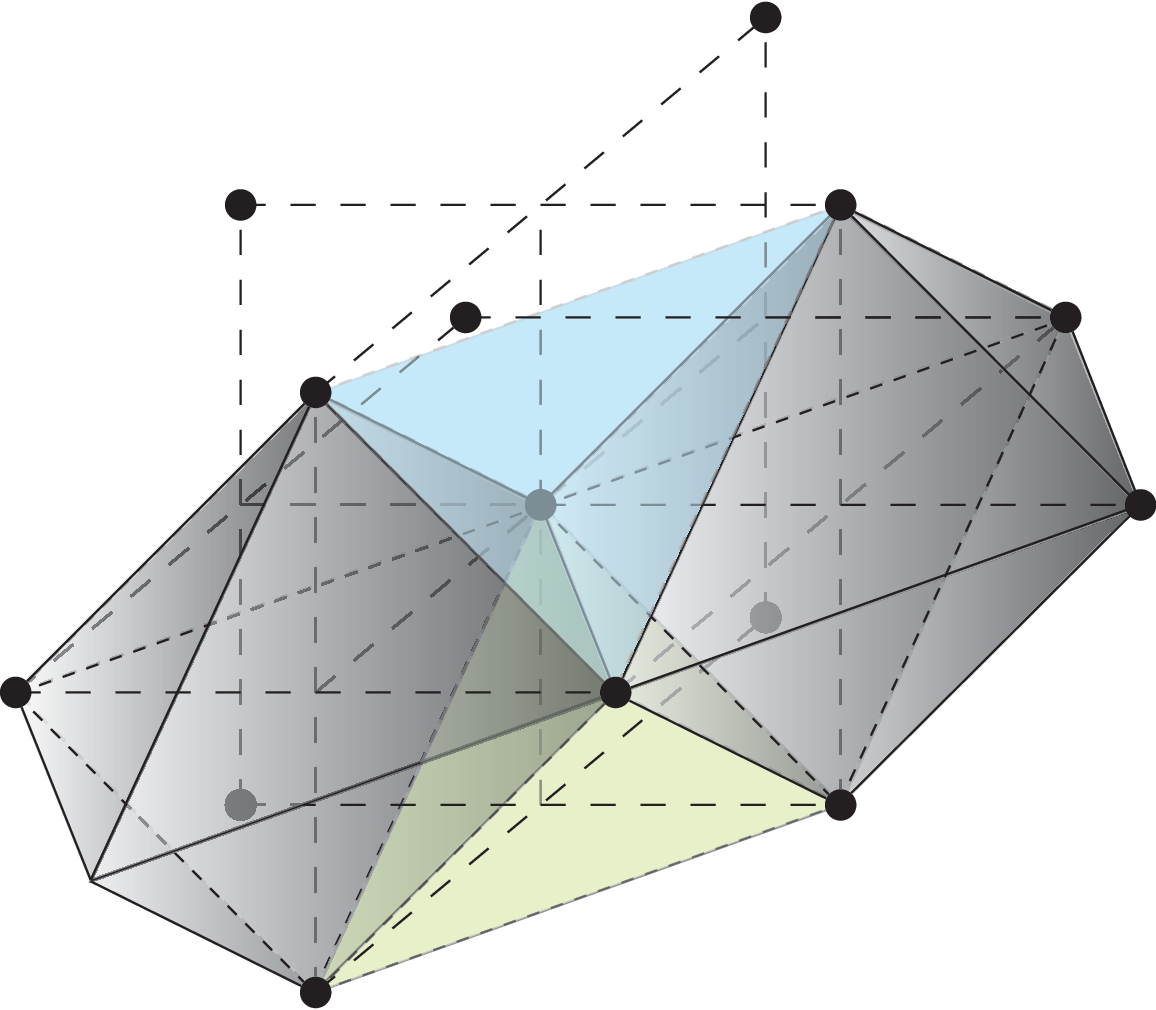}}
    \subfigure[]{\includegraphics[width=0.4\textwidth]{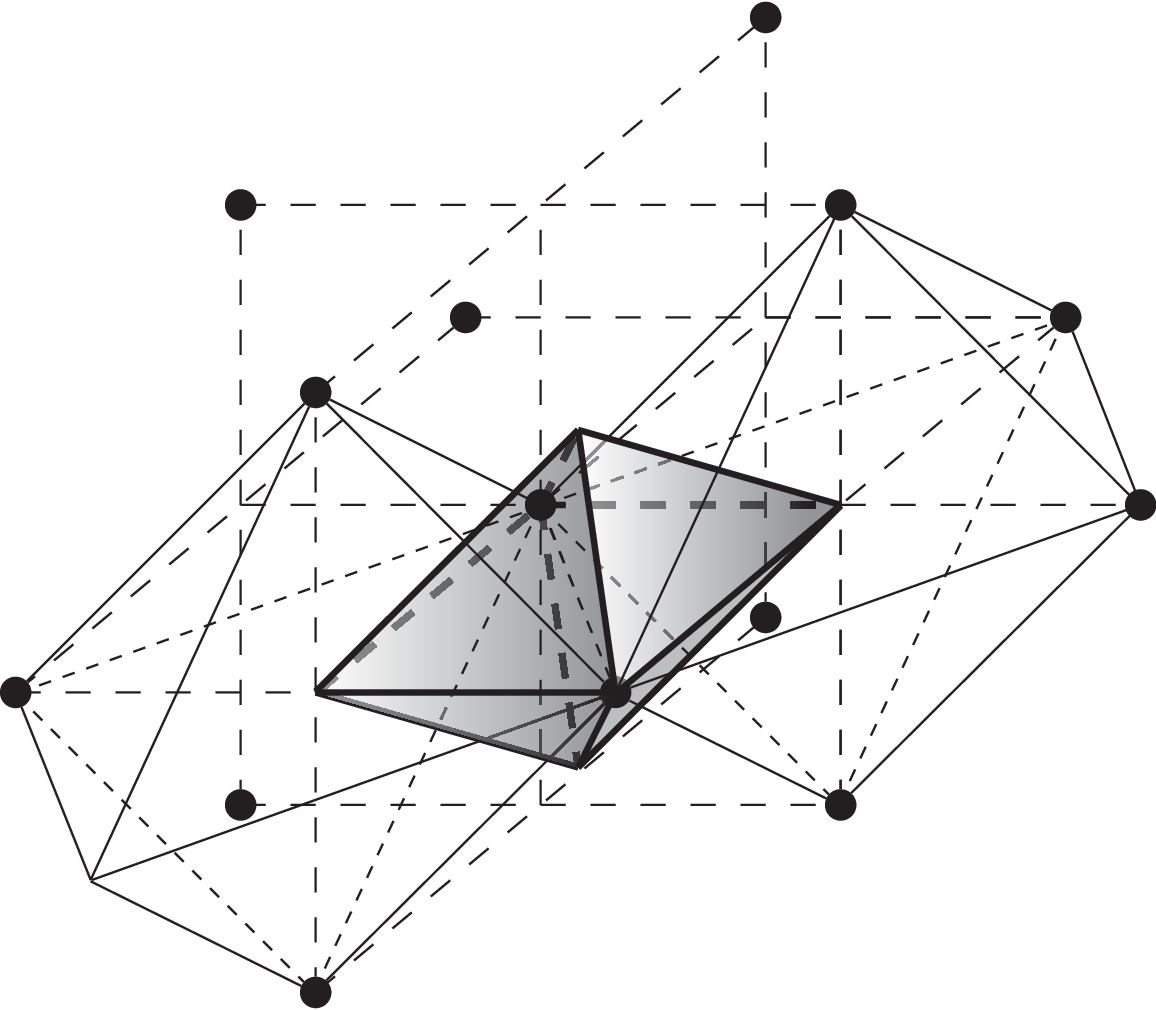}}
    \caption{Octet-truss metamaterial. (a) Space-filling volume cells  \cite([Fig.~8]{Ariza:2005}). (b) Voronoi cell of a bar relative to the volume cells.} \label{S2rFTT}
\end{center}
\end{figure*}

Henceforth, we denote by $E(u;\, \Omega)$ the energy (\ref{w1G0fN}) of the metastructure maximally contained within $\Omega$. For ease of handling, we decompose the energy as
\begin{equation} \label{ifnCt4}
    E(u;\, \Omega) := A(u;\, \Omega) + B(u;\, \Omega) + C(u;\, \Omega) ,
\end{equation}
in terms of an {\sl axial energy}
\begin{equation}
    A(u;\, \Omega)
    =
    \sum_{(m,\beta) \in J_\Omega}
        \frac{{EA}_\beta }{2 L_\beta}
        \Big(
            dv(m,\beta) \cdot d_{1,\beta}
        \Big)^2 ,
\end{equation}
a {\sl bending energy},
\begin{equation} \label{JWMGJR}
\begin{split}
    B(u;\, \Omega)
    & =
    \sum_{(m,\beta) \in J_\Omega}
    \Big\{
        \frac{{GI}_{1,\beta} }{2 L_\beta} \,
        \Big(
            d\theta(m,\beta) \cdot d_{1,\beta}
        \Big)^2 +
\end{split}
\end{equation}
\vglue -0.70truecm
\begin{equation} \nonumber
        \frac{{EI}_{2,\beta} }{2 L_\beta} \,
        \Big(
            d\theta(m,\beta) \cdot d_{2,\beta}
        \Big)^2
        + 
        \frac{{EI}_{3,\beta} }{2 L_\beta} \,
        \Big(
            d\theta(m,\beta) \cdot d_{3,\beta}
        \Big)^2
    \Big\} ,
\end{equation}
and a {\sl coupling energy}
\begin{equation} \label{EhmbE9}
\begin{split}
    &
    C(u;\, \Omega)
    =
    \sum_{(m,\beta) \in J_\Omega}
    \Big\{
    \\ &
        \frac{6 {EI}_{2,\beta}}{L_\beta}
        \Big(
            \frac{dv(m,\beta) \cdot d_{3,\beta}}{L_\beta}
            +
            \bar{\theta}(m,\beta) \cdot d_{2,\beta}
        \Big)^2
        + \\ &
        \frac{6 {EI}_{3,\beta}}{L_\beta}
        \Big(
            \frac{dv(m,\beta) \cdot d_{2,\beta}}{L_\beta}
            -
            \bar{\theta}(m,\beta) \cdot d_{3,\beta}
        \Big)^2
    \Big\} .
\end{split}
\end{equation}

\subsection{Integral representation of the energy} \label{G9O05B}

In order to facilitate the passage to the continuum and make contact with homogenization theory \cite{Cioranescu:1999}
, we rewrite the energy (\ref{ifnCt4}) in integral form. To this end, we introduce a number of auxiliary structures.

Let $P(m,\beta)$ be the {\sl Voronoi cell} of bar $(m,\beta)$ in an infinite metamaterial, i.~e.,
\begin{equation}
\begin{split}
    &
    P(m,\beta)
    =
    \{
        x \in \mathbb{R}^n \, : \,
        {\rm dist}(x,h(m,\beta))
        \leq \\ &
        {\rm dist}(x,h(m',\beta')) ,
        \;\;
        \forall (m',\beta') \in \mathbb{Z}^n\times \{1,\dots, M\}
    \} ,
\end{split}
\end{equation}
i.~e., $P(m,\beta)$ is the locus of points that are closer to be bar $(m,\beta)$ than to any other bar in the metamaterial. The measure
\begin{equation}
    V_\beta := | P(m,\beta) |
\end{equation}
does not depend on $m$ by periodicity and translation invariance. Figs.~\ref{bCb5vv} and \ref{S2rFTT} depict examples of Voronoi cells for a two-dimensional honeycomb structure and a three-dimensional octet-truss, respectively, by way of illustration.

We note that the bar domains $P(m,\beta)$ are polygons, for $n=2$, and polyhedra, for $n=3$, consisting of simplices incident on the spanning segment of the bar, cf.~Figs.~\ref{bCb5vv} and \ref{S2rFTT}. We denote by $\mathcal{T}$ the collection of all such simplices, which together constitute a triangulation of $\mathbb{R}^n$. We note that the vertex set of $\mathcal{T}$ includes---but its larger than---the set of joints of the metamaterial.

Consider now a metastructure maximally contained in $\Omega$. We denote by $\mathcal{T}_\Omega$ the {\sl restriction} of $\mathcal{T}$ to $\Omega$. Precisely, $\mathcal{T}_\Omega$ is the union of the Voronoi cells $P(m,\beta)$ of all bars $(m,\beta) \in J_{\Omega}$. We denote by $[\mathcal{T}_\Omega]$ the domain spanned by $\mathcal{T}_\Omega$. By this construction, $\mathcal{T}_\Omega$ is a triangulation of $[\mathcal{T}_\Omega]$. However, $[\mathcal{T}_\Omega] \neq \Omega$ in general, though the discrepancy is rendered small upon scaling, in a sense to be made precise.

Let $u(l,\alpha)$ be a displacement of the metastructure and ${u}(x)$ a piecewise linear function, not renamed, supported on the triangulation $\mathcal{T}_\Omega$, such that
\begin{equation}\label{eqdefuxula}
    {u}(x(l,\alpha)) = u(l,\alpha) ,
    \quad
    \forall \, (l,\alpha) \in I_{\Omega} .
\end{equation}
Then, $D{u}(x)$ is a piecewise constant function also supported on the triangulation $\mathcal{T}_\Omega$.

We remark that the function $u(x)$ is not uniquely defined by the discrete variables $u(l,\alpha)$, as additional degrees of freedom arise at corners of the Voronoi cells which do not correspond to joints. The value of $u$ at these degrees of freedom can be fixed, for example, by minimizing the mean square oscillation of $Du$ inside each cell of the metastructure, thus rendering the definition unique. In the case of simplicial metamaterials, as for example triangular metamaterials in two dimensions, $Du$ can be made constant on each cell of the metastructure.

Using the identity
\begin{equation} \label{R3aLut}
    L_\beta D{v}(x) \, d_{1,\beta}
    =
    dv(m,\beta) ,
    \qquad
    x \in h(m,\beta) ,
\end{equation}
and the piecewise-constant property of $D{v}(x)$, we obtain
\begin{equation} \label{l5UlMC}
\begin{split}
    &
    A(u;\, \Omega)
    =
    \sum_{(m,\beta) \in J_{\Omega}}
    \\ &
    {
        \frac{1}{V_\beta}
        \int_{P(m,\beta)}
            \frac{{EA}_\beta }{2 L_\beta}
            \Big(
                L_\beta D{v}(x) \, d_{1,\beta} \cdot d_{1,\beta}
            \Big)^2
        \, dx
    } .
\end{split}
\end{equation}
Introducing the piecewise-constant axial moduli
\begin{equation}
    \mathbb{A}(x)
    =
    \sum_{m\in\mathbb{Z}^n} \sum_{\beta=1}^M
        \mathbb{A}_\beta \, \chi_{P(m,\beta)}(x) ,
    \quad
    x \in \mathbb{R}^n ,
\end{equation}
with $\chi_{P(m,\beta)}(x)$ as in (\ref{qe7B0S}) and coefficients
\begin{equation}
    \mathbb{A}_\beta
    =
    \frac{{EA}_\beta L_\beta}{V_\beta} \,
    (d_{1,\beta} \otimes d_{1,\beta}) \otimes (d_{1,\beta} \otimes d_{1,\beta}) ,
\end{equation}
the axial energy (\ref{l5UlMC}) further reduces to the integral form
\begin{equation} \label{FX2Fam}
    A(u;\, \Omega)
    =
    \int_{[\mathcal{T}_\Omega]}
        \frac{1}{2}
        \mathbb{A}(x) \,
        D{v}(x) \cdot D{v}(x)
    \, dx .
\end{equation}
We note that the inhomogeneous moduli $\mathbb{A}(x)$ thus defined are periodic with the periodicity of the metamaterial.

Proceeding likewise, we arrive at the integral form of the bending energy 
\begin{equation} \label{ltriL5}
    B(u;\, \Omega)
    =
    \int_{[\mathcal{T}_\Omega]}
        \frac{1}{2}
        \mathbb{B}(x) \,
        D{\theta}(x) \cdot D{\theta}(x)
    \, dx
\end{equation}
where
\begin{equation}
    \mathbb{B}(x)
    =
    \sum_{m\in\mathbb{Z}^n} \sum_{\beta=1}^M
        \mathbb{B}_\beta \, \chi_{P(m,\beta)}(x) ,
    \quad
    x \in \mathbb{R}^n
\end{equation}
are piecewise-constant bending moduli with coefficients
\begin{equation}
\begin{split}
    \mathbb{B}_\beta
    & =
    \frac{{GI}_{1,\beta} L_\beta}{V_\beta} \,
    (d_{1,\beta} \otimes d_{1,\beta}) \otimes (d_{1,\beta} \otimes d_{1,\beta})
    \\ & +
    \frac{{EI}_{2,\beta} L_\beta}{V_\beta} \,
    (d_{2,\beta} \otimes d_{1,\beta}) \otimes (d_{2,\beta} \otimes d_{1,\beta})
    \\ & +
    \frac{{EI}_{3,\beta} L_\beta}{V_\beta} \,
    (d_{3,\beta} \otimes d_{1,\beta}) \otimes (d_{3,\beta} \otimes d_{1,\beta}) .
\end{split}
\end{equation}

Finally, we turn to the coupling energy. As before, let $u(l,\alpha) = \big( v(l,\alpha);\, \theta(l,\alpha) \big)$ be a displacement of the metastructure. In contrast to the fully conforming interpolation used in the foregoing, we now let ${v}(x)$ be a piecewise linear function and ${\theta}(x)$ a piecewise constant function, not renamed, both supported on the triangulation $\mathcal{T}_\Omega$, and such that
\begin{subequations}
\begin{align}
    &
    v(x(l,\alpha)) = v(l,\alpha) ,
    &
    \forall \, (l,\alpha) \in I_{\Omega} ,
    \\ &
    \theta(x) = \bar{\theta}(m,\beta) ,
    &
    x \in P(m,\beta), \;\; \forall \, (m,\beta) \in J_{\Omega} .
\end{align}
\end{subequations}
As before, $D{v}(x)$ is a piecewise constant function supported on $\mathcal{T}_\Omega$. Using again the identity (\ref{R3aLut}) and the piece\-wise-constant property of $D{v}(x)$, we obtain
\begin{equation} \label{C8Rep0}
    C(u;\, \Omega)
    =
    \sum_{(m,\beta) \in J_\Omega}
    \Big\{ \qquad\qquad\qquad\qquad\qquad\qquad
\end{equation}
\vglue -0.5truecm
\begin{equation} \nonumber
\begin{split}
        &
        \frac{1}{V_\beta}
        \int_{P(m,\beta)}
            \frac{6 {EI}_{2,\beta}}{L_\beta}
            \Big(
                Dv(x) \, d_{1,\beta} \cdot d_{3,\beta}
                +
                \theta(x) \cdot d_{2,\beta}
            \Big)^2
        \, dx
        + \\ &
        \frac{1}{V_\beta}
        \int_{P(m,\beta)}
            \frac{6 {EI}_{3,\beta}}{L_\beta}
            \Big(
                Dv(x) \, d_{1,\beta} \cdot d_{2,\beta}
                -
                \theta(x) \cdot d_{3,\beta}
            \Big)^2
        \, dx
    \Big\} .
\end{split}
\end{equation}
Equivalently, we have
\begin{equation} \label{FTqNWJ}
    C(u;\, \Omega)
    =
    \sum_{(m,\beta) \in J_\Omega}
    \Big\{ \qquad\qquad\qquad\qquad\qquad\qquad
\end{equation}
\vglue -0.5truecm
\begin{equation} \nonumber
\begin{split}
        &
        \frac{1}{V_\beta}
        \int_{P(m,\beta)}
            \frac{6 {EI}_{2,\beta}}{L_\beta}
            \Big(
                \big(
                    Dv(x)
                    -
                    *\theta(x)
                \big)
                \, d_{1,\beta} \cdot d_{3,\beta}
            \Big)^2
        \, dx
        + \\ &
        \frac{1}{V_\beta}
        \int_{P(m,\beta)}
            \frac{6 {EI}_{3,\beta}}{L_\beta}
            \Big(
                \big(
                    Dv(x)
                    -
                    *\theta(x)
                \big)
                \, d_{1,\beta} \cdot d_{2,\beta}
            \Big)^2
        \, dx
    \Big\} ,
\end{split}
\end{equation}
in terms of the Hodge-* operators (\ref{RZBy8m}) and (\ref{dkANTL}). Introducing the piecewise-constant coupling moduli
\begin{equation}
    \mathbb{C}(x)
    =
    \sum_{m\in\mathbb{Z}^n} \sum_{\beta=1}^M
        \mathbb{C}_\beta \, \chi_{P(m,\beta)}(x) ,
    \quad
    x \in \mathbb{R}^n ,
\end{equation}
with coefficients
\begin{equation}
\begin{split}
    \mathbb{C}_\beta
    & =
    \frac{1}{V_\beta}
    \frac{12 {EI}_{2,\beta}}{L_\beta} \,
    (d_{1,\beta} \otimes d_{3,\beta}) \otimes (d_{1,\beta} \otimes d_{3,\beta})
    \\ & +
    \frac{1}{V_\beta}
    \frac{12 {EI}_{3,\beta}}{L_\beta} \,
    (d_{1,\beta} \otimes d_{2,\beta}) \otimes (d_{1,\beta} \otimes d_{2,\beta}) ,
\end{split}
\end{equation}
the coupling energy (\ref{FTqNWJ}) reduces to the integral form
\begin{equation} \label{zKlrnY}
\begin{split}
    &
    C(u;\, \Omega)
    = \\ &
    \int_{[\mathcal{T}_\Omega]}
        \frac{1}{2}
        \mathbb{C}(x) \,
        \big( Dv(x) - *\theta(x) \big)
        \cdot
        \big( Dv(x) - *\theta(x) \big)
    \, dx ,
\end{split}
\end{equation}
which penalizes differences between the deflection-com\-pat\-i\-ble rotations and the mean rotations of the bars. We again note that the inhomogeneous moduli $\mathbb{C}(x)$ thus defined are periodic with the periodicity of the metamaterial.

\begin{figure*}[ht!]
\begin{center}
\includegraphics[width=12.5cm]{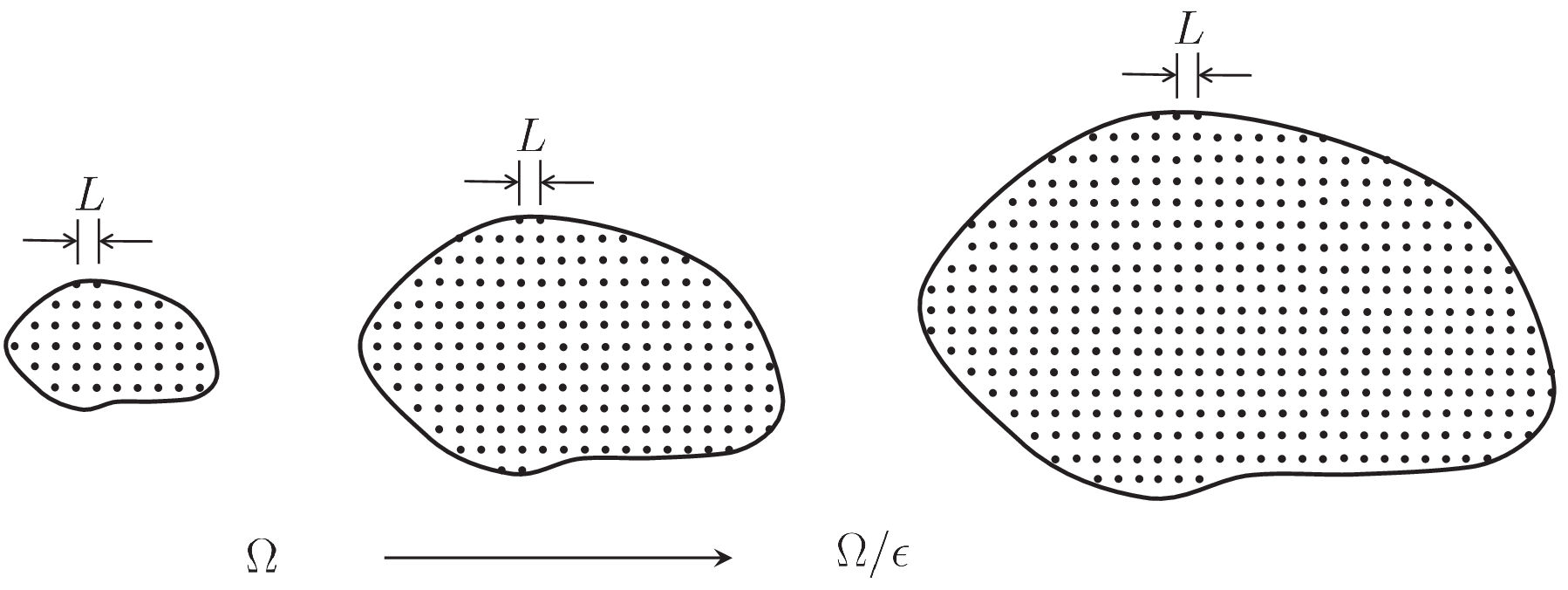}
\end{center}
\caption{Scaling of domain $\Omega$ with constant metamaterial.}
\label{bCb5ny}
\end{figure*}

\subsection{Scaling}
Next, we consider a sequence of metastructures composed of the same metamaterial but spanning increasingly larger self-similar domains $\Omega/\epsilon$, where $\epsilon > 0$ is a small scaling factor, Fig.~\ref{bCb5ny}. The corresponding metastructure is defined by the index set
\begin{equation}
\begin{split}
    &
    J_{\Omega/\epsilon}
    := \\ &
    \{ (m,\beta) \in \mathbb{Z}^n\times \{1,\dots, M\}\, : \,
    h(m,\beta)
    \subset\Omega/\epsilon \} .
\end{split}
\end{equation}
Note that, in this scaling, the geometry of the metamaterial at the unit-cell scale as well and the material properties and the geometry of the cross sections of the bars, remain unchanged along the sequence. In particular, the bar lengths $L_\beta$, cross-sectional areas $A_\beta$ and moments of inertia $I_{i,\beta}$ are {\sl not} scaled according to $\epsilon$. By contrast, the number of bars in the metastructure grows as $\# I_{\Omega/\epsilon} \sim O(\epsilon^{-n})$ and the metastructure spans an infinite metamaterial in the {\sl infinite body limit} $\epsilon \to 0$.

With $\epsilon \to 0$, we consider the sequence of energy functions
\begin{equation} \label{3R5GsZ}
    E_\epsilon(u;\, \Omega)
    =
    \epsilon^{n} \, E(u;\, \Omega/\epsilon) ,
\end{equation}
with
\begin{equation} \label{32yGNG}
    x_\epsilon =  \epsilon^{-1} \, x,
    \quad
    v_\epsilon(x_\epsilon)
    =
    \epsilon^{-1} \, v(x) ,
    \quad
    \theta_\epsilon(x_\epsilon)
    = \theta(x) .
\end{equation}
The factor $\epsilon^{n}$ in (\ref{3R5GsZ}) accounts for the expected elasticity scaling of the limiting energy and is included to ensure a proper limit.

\begin{figure*}[ht!]
\begin{center}
    \includegraphics[width=12.5cm]{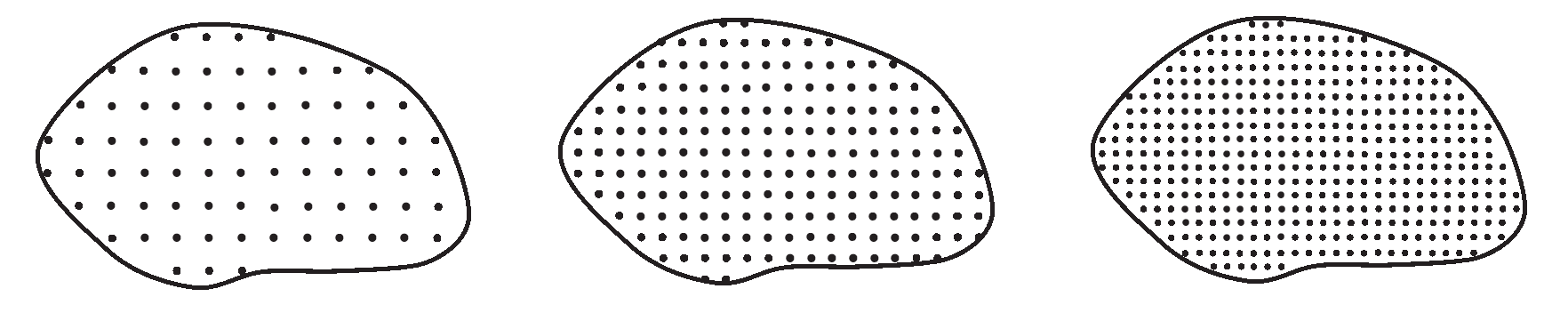}
\end{center}
\caption{Scaling of metamaterial with constant domain $\Omega$.}
\label{bCb5n2}
\end{figure*}

From (\ref{FX2Fam}), (\ref{ltriL5}) and (\ref{zKlrnY}), we have the identity
\begin{equation} \nonumber
    E(u;\, \Omega/\epsilon)
    =
    A(u;\, \Omega/\epsilon)
    +
    B(u;\, \Omega/\epsilon)
    +
    C(u;\, \Omega/\epsilon)
    =
\end{equation}
\vglue -0.85truecm
\begin{equation} \label{6kWfZ1}
\begin{split}
    &
    \int_{[\mathcal{T}_{\Omega/\epsilon}]}
        \frac{1}{2}
        \mathbb{A}(x) \,
        D{v}(x) \cdot D{v}(x)
    \, dx
    + \qquad\qquad {} \;\;\;\;\;\;\;\; \\ &
    \int_{[\mathcal{T}_{\Omega/\epsilon}]}
        \frac{1}{2}
        \mathbb{B}(x) \,
        D{\theta}(x) \cdot D{\theta}(x)
    \, dx
    + \qquad\qquad {} \;\;\;\;\;\;\;\;
\end{split}
\end{equation}
\vglue -0.625truecm
\begin{equation} \nonumber
    \int_{[\mathcal{T}_{\Omega/\epsilon}]}
        \frac{1}{2}
        \mathbb{C}(x) \,
        \big( Dv(x) - *\theta(x) \big)
        \cdot
        \big( Dv(x) - *\theta(x) \big)
    \, dx ,
\end{equation}
for displacement fields $u(x)$ in the class of piecewise-linear or piecewise-constant functions supported on $\mathcal{T}$ defined in Section~\ref{G9O05B}.

The structure of (\ref{6kWfZ1}) is noteworthy. Thus, the energy consists of three terms: an axial energy; a bending energy, including torsion; and a coupling energy that couples deflection-compatible rotations and free rotations. Energies of the form (\ref{6kWfZ1}) were termed {\sl micropolar} by \cite{Eringen:1966} and set forth an example of a theory of {\sl generalized continua}. Such theories are often used to introduce an {\sl ad hoc} micromechanical length into the material law that renders the energy {\sl non-local}, in an effort to account for size effects. By contrast, in the present context the micropolar terms are physical and introduced into the energy by bending.

Rescaling the domain back to $\Omega$ according to (\ref{3R5GsZ}) and (\ref{32yGNG}), we obtain
\begin{equation} \label{fMCF9Q}
    E_\epsilon(u;\, \Omega)
    =
\qquad\qquad\qquad\qquad\qquad\qquad\qquad\qquad\quad
\end{equation}
\vglue -0.85truecm
\begin{equation} \nonumber
\begin{split}
    &
    \int_{\epsilon [\mathcal{T}_{\Omega/\epsilon}]}
        \frac{1}{2}
        \mathbb{A}(x/\epsilon) \,
        D{v}(x) \cdot D{v}(x)
    \, dx
    + \\ &
    \epsilon^2 \,
    \int_{\epsilon [\mathcal{T}_{\Omega/\epsilon}]}
        \frac{1}{2}
        \mathbb{B}(x/\epsilon) \,
        D{\theta}(x) \cdot D{\theta}(x)
    \, dx
    + \\ &
    \int_{\epsilon [\mathcal{T}_{\Omega/\epsilon}]}
        \frac{1}{2}
        \mathbb{C}(x/\epsilon) \,
        \big( Dv(x) - *\theta(x) \big)
        \cdot
        \big( Dv(x) - *\theta(x) \big)
    \, dx ,
\end{split}
\end{equation}
where we have dropped the prime for simplicity and the new displacement fields $u(x)$ are in the class of piecewise-linear or piecewise-constant functions supported on the scaled triangulation $\epsilon \, \mathcal{T}$.

\subsection{Continuum limit}

We wish to ascertain the asymptotic behavior of the {\sl forced} energy-minimizing configurations $(u_\epsilon)$ of the discrete energies (\ref{fMCF9Q}) as $\epsilon \to 0$, or {\sl continuum limit}. Specifically, we wish to determine whether a continuum limiting energy $E_0$ exists with the property that the sequence of minimum discrete potential energies converges to the minimum continuum potential energy {\sl for all admissible forcings}. Then, we say that $E_0$ is the variational limit, or $\Gamma$-limit, of the sequence $(E_\epsilon)$ \cite{dalmaso:1993}).

Thus, suppose that the metastructure is acted upon by macroscopic distributed actions defined by the functions
\begin{equation}
    f_0(x) := \big( q_0(x) ;\, m_0(x) \big),
\end{equation}
where $q_0(x)$ and $m_0(x)$ are distributed forces and moments, or torques, per unit volume, respectively. The work of the applied actions is then
\begin{equation} \label{XVOWv4}
    \langle f_0,\, u \rangle
    :=
    \frac{V}{N} \,
    \sum_{l\in\mathbb{Z}^n}
        \sum_{\alpha=1}^N
        f_0(x(l,\alpha)) \cdot u(l,\alpha) ,
\end{equation}
where $V$ is the volume of the unit cell of the spanning Bravais lattice of the metamaterial and the factor $V/N$ represents the volume per joint class. Its value follows from the property that all joint classes have the same number of joints per unit volume, or joint density.

Consider now a sequence $(E_\epsilon)$ of scaled energies defined in (\ref{fMCF9Q}). Assume now that the sequence of energies is {\sl uniformly stable}, or {\sl equicoercive}. Specifically, we assume that there is a constant $C > 0$ such that
\begin{equation} \label{9JbQjQ}
    E_\epsilon(u) \geq C \, \| u \|^2 ,
\end{equation}
for some suitable {\sl energy norm} $\|\cdot\|$. This equicoercivity condition expresses the property that $E_\epsilon(u)$ is bounded below for all $\epsilon$ by a {\sl structurally stable} quadratic energy, i.~e., a positive quadratic energy devoid of zero-energy modes. This property effectively prevents the emergence of finely oscillating unstable modes that corrupt the solution, at no cost in energy, as $\epsilon\to 0$.

The norm entering~\eqref{9JbQjQ} needs to be chosen carefully, and a suitable choice may not be immediately apparent for some metastructures. Besides the usual fact that it should only control $u$ modulo rigid-body motions, one needs to consider several conditions. First, the metamaterial needs to be appropriately truncated at the boundary, eliminating for example dangling bonds from the domain on which $\|\cdot\|$ is computed. Second, the norm needs to be strong enough to ensure continuity of the energy up to boundary terms. Third, the ancillary degrees of freedom need to be eliminated as discussed after~\eqref{eqdefuxula}, and the resulting metamaterial needs to be sufficiently rigid to provide coercivity. For example, for a triangular metamaterial in two dimensions, the axial energy immediately provides control of the $L^2$ norm of the strain $(Dv+Dv^T)/2$, which leads to control of $v$ and $Dv$ via Korn's inequality and appropriate boundary data. Using this property, the coupling energy gives control of $\theta$, and finally the bending energy gives control of $D\theta/\epsilon$. The situation for the honeycomb lattice is more complex, as mechanisms exist for the axial energy, which correspond to zero-axial-energy deformations in addition to rigid-body motions. These zero-energy modes need to be stabilized by the bending and coupling energies in the long-wavelength limit relevant to the macroscopic material behavior.

With these provisos, the equilibrium configuration $u_\epsilon$ of the metastructure then follows by minimizing the potential energy
\begin{equation} \label{cLMgtE}
    F_\epsilon(u)
    :=
    E_\epsilon(u)
    -
    \langle f_0,\, u \rangle .
\end{equation}
Denote by $m(E_\epsilon; f_0)$ the corresponding minimum potential energy. Then, the variational limit, or $\Gamma$-limit, $E_0$ of the sequence $(E_\epsilon)$ is characterized by the property that
\begin{equation} \label{0LG0G2}
    \lim_{\epsilon\to 0}
    m(E_\epsilon; f_0)
    =
    m(E_0; f_0) ,
\end{equation}
for all admissible forcings $f_0(x)$, where $m(E_0; f_0)$ is the minimum of the limiting potential energy
\begin{equation} \label{aB9UpY}
    F_0(u_0) := E_0(u_0) - \langle f_0,\, u_0 \rangle ,
\end{equation}
(cf., e.~g., \cite[Theorem 13.5, Corollary 13.7]{dalmaso:1993}). Variational convergence, together with the equicoercivity condition (\ref{9JbQjQ}), in turn ensures the weak convergence of the sequence of minimizers $(u_\epsilon)$ to the minimizer $u_0$ of $E_0$ for all admissible forcings $f_0$.

In view of the structure of the energy (\ref{ifnCt4}), the corresponding variational limit as $\epsilon\to 0$ extends variational {\sl discrete-to-continuum problems} \cite{BraidesGelli2006,Alicandro:2004} to account for bending. Alternatively, appealing to the density of piecewise linear and constant displacements over the increasingly finer triangulations $\epsilon \mathcal{T}$, the energy (\ref{fMCF9Q}) can be extended (by infinity) to continuum displacement fields without altering the continuum limit. The resulting limiting problem then falls within the class of variational {\sl homogenization problems} (cf., e.~g., \cite{Cioranescu:1999}).

In addition, in view of the $\epsilon^2$ factor in (\ref{fMCF9Q}) we expect the bending energy to be negligible in the limit, and the limiting continuum energy to be of the form
\begin{equation} \label{lxJTdK}
    E_0(u_0;\, \Omega)
    =
    \int_{\Omega}
        W_0\big(Dv_0(x),\, \theta_0(x) \big)
    \, dx ,
\end{equation}
for some effective energy density $W_0( Dv_0,\, \theta_0 )$, to be determined, and for a certain class of sufficiently regular and finite domains $\Omega$. In particular, the effective energy density $W_0( Dv_0,\, \theta_0 )$ is expected to be independent of $\Omega$ within the admissible class.

\subsection{Discrete Fourier transform analysis}

Since, by assumption, the limiting energy density is independent of the domain $\Omega$, we can conveniently identify $W_0( Dv_0,\, \theta_0 )$ by considering metamaterials of infinite extent. In order to make calculations explicit, it proves convenient to revert to the discrete form (\ref{w1G0fN}) or the energy, with $J = \mathbb{Z}^n \times \{1,\cdots,M\}$.

We can exploit the periodicity of the metamaterial by testing the energy with discrete harmonic displacements
\begin{equation}
    u(l,\alpha)
    =
    \frac{1}{(2\pi)^n}
    \int_B
        \hat{u}(k,\alpha) \, {\rm e}^{i k \cdot x(l,\alpha)}
    \, dk ,
\end{equation}
where $B$ is the Brillouin zone of the metastructure, $x(l,\alpha)$ are the joint coordinates (\ref{M6JjBD}) and $\hat{u}(k,\alpha)$ is the discrete Fourier transform of the discrete displacement field $u(l,\alpha)$, cf.~\ref{TwWqmj}. A straightforward calculation using Parseval's theorem (\ref{pB0rWJ}) then gives
\begin{equation} \label{59joLe}
\begin{split}
    &
    A(u)
    = \\ &
    \sum_{\beta=1}^M
    \frac{1}{(2\pi)^n}
    \int_B
        \frac{1}{V}
        \Big\{
            \frac{{EA}_\beta }{2 L_\beta}
            \big|
                d\hat{v}(k,\beta) \cdot d_{1,\beta}
            \big|^2
        \Big\}
    \, dk ,
\end{split}
\end{equation}
for the axial energy,
\begin{equation} \label{XLBfZ6}
\begin{split}
    &
    B(u)
    = 
    \sum_{\beta=1}^M
    \frac{1}{(2\pi)^n}
    \int_B
        \frac{1}{V}
        \Big\{ \\ & \qquad
            \frac{{GI}_{1,\beta} }{2 L_\beta} \,
            \big|
                d\hat{\theta}(k,\beta) \cdot d_{1,\beta}
            \big|^2
            + \\ & \qquad
            \frac{{EI}_{2,\beta} }{2 L_\beta} \,
            \big|
                d\hat{\theta}(k,\beta) \cdot d_{2,\beta}
            \big|^2
            + \\ & \qquad
            \frac{{EI}_{3,\beta} }{2 L_\beta} \,
            \big|
                d\hat{\theta}(k,\beta) \cdot d_{3,\beta}
            \big|^2
        \Big\}
    \, dk ,
\end{split}
\end{equation}
for the bending energy, and
\begin{equation} \label{OVv27j}
\begin{split}
    &
    C(u)
    =
    \sum_{\beta=1}^M
    \frac{1}{(2\pi)^n}
    \int_B
        \frac{1}{V}
        \Big\{ \\ &
            \frac{6 {EI}_{2,\beta}}{L_\beta}
            \Big|
                \frac{d\hat{v}(k,\beta) \cdot d_{3,\beta}}{L_\beta}
                +
                \bar{\hat{\theta}}(k,\beta) \cdot d_{2,\beta}
            \Big|^2
            + \\ &
            \frac{6 {EI}_{3,\beta}}{L_\beta}
            \Big|
                \frac{d\hat{\theta}(k,\beta) \cdot d_{2,\beta}}{L_\beta}
                -
                \bar{\hat{\theta}}(k,\beta) \cdot d_{3,\beta}
            \Big|^2
        \Big\}
    \, dk ,
\end{split}
\end{equation}
for the coupling energy, where $V$ is the volume of the unit cell of the spanning Bravais lattice of the metamaterial and we write
\begin{subequations} \label{gy0mqM}
\begin{align}
    &
    d\hat{v}(k,\beta)
    =
    {\rm e}^{\frac{i}{2} k \cdot dx_\beta}
    \hat{v}(k,\beta^+)
    -
    {\rm e}^{- \frac{i}{2} k \cdot dx_\beta}
    \hat{v}(k,\beta^-) ,
    \\ &
    d\hat{\theta}(k,\beta)
    =
    {\rm e}^{\frac{i}{2} k \cdot dx_\beta}
    \hat{\theta}(k,\beta^+)
    -
    {\rm e}^{- \frac{i}{2}  k \cdot dx_\beta}
    \hat{\theta}(k,\beta^-) ,
    \\ &
    \bar{\hat{\theta}}(k,\beta)
    =
    \frac{1}{2}
    \Big(
        {\rm e}^{\frac{i}{2}  k \cdot dx_\beta}
        \hat{\theta}(k,\beta^+)
        +
        {\rm e}^{- \frac{i}{2}  k \cdot dx_\beta}
        \hat{\theta}(k,\beta^-)
    \Big) .
\end{align}
\end{subequations}
Collecting all joint displacements into an array, not renamed,
\begin{equation} \label{qbNEqL1}
    \hat{u}(k) := \big( \hat{u}(k,1), \dots, \hat{u}(k,N) \big) ,
\end{equation}
the energy then takes the compact form
\begin{equation} \label{oEf3Zx}
    E(u)
    =
    \frac{1}{(2\pi)^n}
    \int_B
        \frac{1}{2}
        \mathbb{D}(k)
        \, \hat{u}(k) \cdot \hat{u}^*(k)
    \, dk ,
\end{equation}
in terms of the {\sl dynamical matrix} $\mathbb{D}(k)$ of the metamaterial.

We note that $\mathbb{D}(k)$ is Hermitian, positive definite and depends smoothly on the wavevector $k$. We additionally note that $\mathbb{D}(k)$ follows {\sl explicitly} from (\ref{59joLe}), (\ref{XLBfZ6}), (\ref{OVv27j}) and (\ref{gy0mqM}), by a simple rearrangement of terms, as a function of the beam properties: $EA_\beta$, $GI_{1,\beta}$, $EI_{2,\beta}$ and $EI_{3,\beta}$; and the metamaterial geometry: $L_\beta$, $dx_\beta$, $d_{1,\beta}$, $d_{2,\beta}$ and $d_{3,\beta}$.

In the Fourier representation, an application of Parseval's theorem (\ref{pB0rWJ}) to the work function (\ref{XVOWv4}) further gives
\begin{equation} \label{W8oi93}
    \langle f_0,\, u \rangle
    :=
    \frac{1}{(2\pi)^n}
    \int_B
        \hat{f}_0(k)
        \Big(
            \frac{1}{N}
            \sum_{\alpha=1}^N u^*(k,\alpha)
        \Big)
    \, dk ,
\end{equation}
where $\hat{f}_0(k)$ is the ordinary Fourier transform of $f_0(x)$.

The scaling (\ref{9lQ7Sk}) defines the sequence of scaled energies
\begin{equation} \label{tt54ez}
    E_\epsilon(u)
    =
    \frac{1}{(2\pi)^n}
    \int_{B/\epsilon}
        \frac{1}{2}
        \mathbb{D}_\epsilon(k)
        \, \hat{u}(k) \cdot \hat{u}^*(k)
    \, dk ,
\end{equation}
where the scaled dynamical matrix $\mathbb{D}_\epsilon(k)$ is defined by the property
\begin{equation} \label{3ll62m}
\begin{split}
  &
  \mathbb{D}_\epsilon(k) \, \zeta \cdot \zeta^*
  =
  \epsilon^{-2} \,
  \mathbb{D}(\epsilon k) \,
    (\xi; \, \epsilon \eta)
    \cdot
    (\xi^*; \, \epsilon \eta^*) ,
  \\ &
  \forall \zeta := (\xi;\,\eta) \in [\mathbb{C}^n]^N \times [\mathbb{C}^{n(n-1)/2}]^N ,
  \quad
  k \in B/\epsilon ,
\end{split}
\end{equation}
Likewise, as $\epsilon\to 0$ the work function (\ref{W8oi93}) can be expressed independently of $\epsilon$ as
\begin{equation}
    \langle f_0,\, u \rangle
    :=
    \frac{1}{(2\pi)^n}
    \int_{\mathbb{R^n}}
        \hat{f}_0(k)
        \Big(
            \frac{1}{N}
            \sum_{\alpha=1}^N u^*(k,\alpha)
        \Big)
    \, dk .
\end{equation}
provided that $\hat{u}(k)$ is extended to zero outside $B/\epsilon$ and
$\hat{f}_0(k)$ is kept fixed throughout the sequence.

Suppose now that the metamaterial is loaded by macroscopic forces $f_0(x)$. Then, from the Fourier representations (\ref{tt54ez}) and (\ref{W8oi93}) the corresponding equilibrium equations are
\begin{equation} \label{ZqVspb}
    \mathbb{D}_\epsilon(k) \, \hat{u}_\epsilon(k)
    =
    \hat{f}(k) ,
    \quad
    k \in B/\epsilon ,
\end{equation}
where again we collect all the joint forces into an array
\begin{equation} \label{qbNEqLn2}
    \hat{f}(k)
    :=
    \Big(
        \frac{1}{N}
        \hat{f}_0(k),
        \stackrel{N}{\dots},
        \frac{1}{N}
        \hat{f}_0(k)
    \Big)
    :=
    L \hat{f}_0(k) .
\end{equation}
We note that the operator $L$ thus defined {\sl localizes} the continuum forces to each of the joint classes.

In this representation, an appropriate form of the equicoercivity condition (\ref{9JbQjQ}) is that
\begin{equation} \label{P5wyEw}
    \mathbb{D}_\epsilon(k)
    \, \zeta \cdot \zeta^*
    \geq
    C ( |k|^2 \, |\xi|^2 + | \eta |^2 ) ,
\end{equation}
for all $\zeta := (\xi;\,\eta) \in \mathbb{C}^n \times \mathbb{C}^{n(n-1)/2}$ and some $C > 0$. We note that this condition effectively implies a corresponding energy norm $\|\cdot\|$, eq.~(\ref{9JbQjQ}). Under these conditions, the matrix $\mathbb{D}_\epsilon(k)$ is non-singular and the equilibrium problem (\ref{ZqVspb}) can be solved pointwise for $k \in B/\epsilon$, with the result
\begin{equation} \label{jQ2xKs}
    \hat{u}_\epsilon(k)
    =
    \mathbb{D}_\epsilon^{-1}(k) \, \hat{f}(k) ,
    \quad
    k \in B/\epsilon ,
\end{equation}
which {\sl explicitly} characterizes the sequence of energy-min\-i\-mizing displacements $(\hat{u}_\epsilon)$ over their corresponding Brillouin zones $B/\epsilon$.

We note that we can equivalently extend the functions $\hat{u}_\epsilon(k)$ to all of $\mathbb{R}^n$ by zero. The corresponding inverse Fourier transforms $u_\epsilon(x)$ are then said to be {\sl band-limited} and are referred to as the Whittaker-Shannon interpolation of the discrete displacements $u_\epsilon(l,\alpha)$ \cite{Mallat:2009, Espanol:2013}.

The corresponding minimum potential energy is
\begin{equation}
    m(E_\epsilon; f_0)
    :=
    -
    \frac{1}{(2\pi)^n}
    \int_{B/\epsilon}
        \frac{1}{2}
        \mathbb{D}_\epsilon^{-1}(k)
        \, \hat{f}(k) \cdot \hat{f}^*(k)
    \, dk .
\end{equation}
Passing to the limit as in (\ref{0LG0G2}) with the aid of (\ref{P5wyEw}) and dominated convergence, we obtain
\begin{equation}
\begin{split}
    &
    \lim_{\epsilon\to 0}
    m(E_\epsilon; f_0)
    :=
    m(E_0; f_0)
    =
    - E_0^*(f_0)
    = \\ &
    -
    \frac{1}{(2\pi)^n}
    \int_{\mathbb{R}^n}
        \frac{1}{2}
        \mathbb{D}_0^{-1}(k)
        \, \hat{f}_0(k) \cdot \hat{f}_0^*(k)
    \, dk ,
\end{split}
\end{equation}
where
\begin{equation} \label{U2E0YW}
    \mathbb{D}_0(k)
    =
    \Big(
        \lim_{\epsilon\to 0}
        L^T \mathbb{D}_\epsilon^{-1}(k) L
    \Big)^{-1} ,
    \quad
    k \in \mathbb{R}^n ,
\end{equation}
is the effective dynamical matrix in the continuum limit and $E_0^*(f_0)$ is the dual continuum energy. A final Legendre transform gives the effective continuum energy as
\begin{equation} \label{Lc9K2l}
    E_0(u_0)
    =
    \frac{1}{(2\pi)^n}
    \int_{\mathbb{R}^n}
        \frac{1}{2}
        \mathbb{D}_0(k)
        \, \hat{u}_0(k) \cdot \hat{u}_0^*(k)
    \, dk ,
\end{equation}
for all continuum displacements $u_0(x)$.

The limiting equilibrium problem for the continuum displacements $u_0(x)$ consists of minimizing the homogenized potential energy (\ref{aB9UpY}) for given forcing $f_0(x)$. In the Fourier representation, the solution is
\begin{equation}
    \hat{u}_0(k) = \mathbb{D}_0^{-1}(k) \, \hat{f}_0(k),
    \quad
    k \in \mathbb{R}^n ,
\end{equation}
which combined with (\ref{jQ2xKs}) gives
\begin{equation}
    \hat{u}_\epsilon(k)
    =
    \mathbb{D}_\epsilon^{-1}(k) L \mathbb{D}_0(k) \hat{u}_0(k),
    \quad
    k \in B/\epsilon .
\end{equation}
This relation shows how the continuum displacements $u_0(x)$ are approximated by discrete displacements $u_\epsilon(l,\alpha)$ along the sequence which are the result of energy relaxation at the unit-cell level.

We additionally note that, for metamaterials with one class of joints only, $N=1$, the matrix $L$ reduces to the identity, whereupon (\ref{U2E0YW}) in turn reduces to \cite{Ariza:2005}
\begin{equation} \label{byu9Ue}
    \mathbb{D}_0(k)
    =
    \lim_{\epsilon\to 0}
    \mathbb{D}_\epsilon(k) ,
\end{equation}
i.~e., $\mathbb{D}_0(k)$ is the pointwise limit of $\mathbb{D}_\epsilon(k)$ as $\epsilon\to 0$, or {\sl long wavelength} limit.

\subsection{Structure of the continuum limit}

We proceed to elucidate the structure of the limiting continuum energy (\ref{Lc9K2l}) as a quadratic form. A lengthy but straightforward calculation, consigned to \ref{KA4bct} for the sake of continuity, shows that
\begin{equation} \label{qwzrN3}
\begin{split}
    &
    \mathbb{D}_0(\lambda k)
    \, {\zeta}_0 \cdot {\zeta}_0^*
    = 
    \mathbb{D}_0(k) \,
    (\lambda \xi_0;\, \eta_0)
    \cdot
    (\lambda \xi_0^*;\, \eta_0^*) .
\end{split}
\end{equation}
for every $\lambda > 0$ and
\begin{equation} \label{8t4tSW}
    \zeta_0 := (\xi_0; \, \eta_0)
    \in \mathbb{C}^n \times \mathbb{C}^{n(n-1)/2} ,
\end{equation}
i.~e., rescaling $k$ by $\lambda$ in the continuum dynamical matrix $\mathbb{D}_0(k)$ is equivalent to rescaling the deflection amplitude $\xi_0$ by $\lambda$, while simultaneously keeping the rotation amplitude $\eta_0$ invariant.

From this property, it follows that the energy density admits the representation
\begin{equation}
    \frac{1}{2}
    \mathbb{D}_0(k)
    \, {\zeta}_0 \cdot {\zeta}_0^*
    =
    W_0( ik\, \xi_0,\, \eta_0) ,
\end{equation}
where $W_0$ is a quadratic form, and the energy (\ref{Lc9K2l}) becomes
\begin{equation}
    E_0(u_0)
    =
    \frac{1}{(2\pi)^n}
    \int_{\mathbb{R}^n}
        W_0\big( ik\, \hat{v}_0(k),\, \hat{\theta}_0(k) \big)
    \, dk .
\end{equation}
An inverse Fourier transform then gives the energy as
\begin{equation} \label{J3gwXW}
    E_0(u_0)
    =
    \int_{\mathbb{R}^n}
        W_0\big(Dv_0(x),\, \theta_0(x) \big)
    \, dx ,
\end{equation}
or (\ref{lxJTdK}) for a domain $\Omega$, as surmised.

Further restrictions on $W_0\big(Dv_0(x),\, \theta_0(x) \big)$ result from material-frame invariance, eq.~(\ref{UdCT8i}), which requires that
\begin{equation} \label{VHFvrd}
    W_0\big(\beta,\, \theta \big)
    =
    W_0\big(\beta + w,\, \theta + *w\big) ,
\end{equation}
for all deformation gradients $\beta \in \mathbb{R}^{n\times n}$, rotation angles $\theta \in \mathbb{R}^{n(n-1)/2}$ and rotation matrices $w \in so(n)$. Choosing $w = -{\rm skw} \beta$, eq.~(\ref{VHFvrd}) specializes to
\begin{equation} \label{sNoddM}
    W_0\big(\beta,\, \theta \big)
    =
    W_0\big(\varepsilon ,\, \theta - *{\rm skw} \beta \big)
\end{equation}
for some other function $W_0$, not renamed, of the strains $\varepsilon = {\rm sym} \beta$ and the difference between the free rotations $\theta$ and the deflection-compatible rotations $*{\rm skw} \beta$. Inserting (\ref{sNoddM}) into (\ref{J3gwXW}), we obtain the representation
\begin{equation} \label{J3gwXWb}
    E_0(u_0)
    =
    \int_{\Omega}
        W_0\big(\varepsilon_0(x),\, \theta_0(x) - *{\rm skw} Dv_0(x)\big)
    \, dx ,
\end{equation}
where
\begin{equation}
    \varepsilon_0(x) = {\rm sym} Dv_0(x)
\end{equation}
are the local strains.

We thus conclude that, to lowest order, the effective continuum energy (\ref{J3gwXWb}) of linear metastructures is a special case of linear {\sl micropolar elasticity}, in the sense of  \cite{Eringen:1966}, in which the energy density is independent of the curvature, or bending strain, $D\theta_0(x)$. Because of this special form, the energy (\ref{J3gwXW}) is {\sl local}, i.~e., it exhibits the local elasticity scaling
\begin{equation}
    E_0(u_0';\; \Omega')
    =
    \lambda^n
    E_0(u_0;\; \Omega) ,
\end{equation}
under the transformation
\begin{equation}
    x' = \lambda x ,
    \quad
    v_0' = \lambda v_0 ,
    \quad
    \theta_0' = \theta_0 .
\end{equation}
In particular, to lowest-order the continuum energy does not account for {\sl size effects}.

We additionally note that the variational convergence properties of the energy carry over to potential energies involving forcing terms. We also expect the energy-min\-i\-mizing configurations of the continuum problem to restrict properly to the boundary in the sense of traces, provided that the metastructure is terminated at the boundary in an appropriate way. Under these conditions, the continuum energy-minimizing configurations satisfy the boundary value problem
\begin{subequations}
\begin{align}
    &
    \operatorname{\rm div}
    \partial_{Dv_0} W_0\big(Dv_0(x),\, \theta_0(x)\big)
    +
    q_0(x)
    =
    0 ,
    \;\;
    x \in \Omega ,
    \\ & \label{Y7HLfb}
    \partial_{\theta_0} W_0\big(Dv_0(x),\, \theta_0(x) \big)
    =
    m_0(x) ,
    \quad
    x \in \Omega ,
    \\ &
    u_0(x) = g_0(x) ,
    \quad
    x \in \Gamma_D ,
    \\ &
    \partial_{Dv_0} W_0\big( Dv_0(x),\, \theta_0(x) \big) \, \nu = h_0(x) ,
    \quad
    x \in \Gamma_N ,
\end{align}
\end{subequations}
where $q_0(x)$, $m_0(x)$, $g_0(x)$ and $h_0(x)$ are continuum body forces, body moments, prescribed deflections and prescribed tractions, respectively; $\Gamma_D$ and $\Gamma_N$ are the displacement, or {\sl Dirichlet}, and the traction, or {\sl Neumann}, parts of the boundary, respectively; and $\nu$ is the outward unit normal at the boundary.

\section{Examples of metamaterials} \label{vX3HxB}

We illustrate the properties of the continuum limit by means of elementary examples. We further illustrate the range and scope of the theory by means of two common examples of  metamaterials: the two-dimensional honeycomb lattice; and the three-dimensional octet-truss. A detailed {\tt Mathematica} (Wolfram Research, Inc.) implementation of the examples can be found in the supplementary materials.

\subsection{One-dimensional metamaterial}\label{TBhq2c}

\begin{figure}[h]
 \begin{center}
 \includegraphics[width=0.35\textwidth]{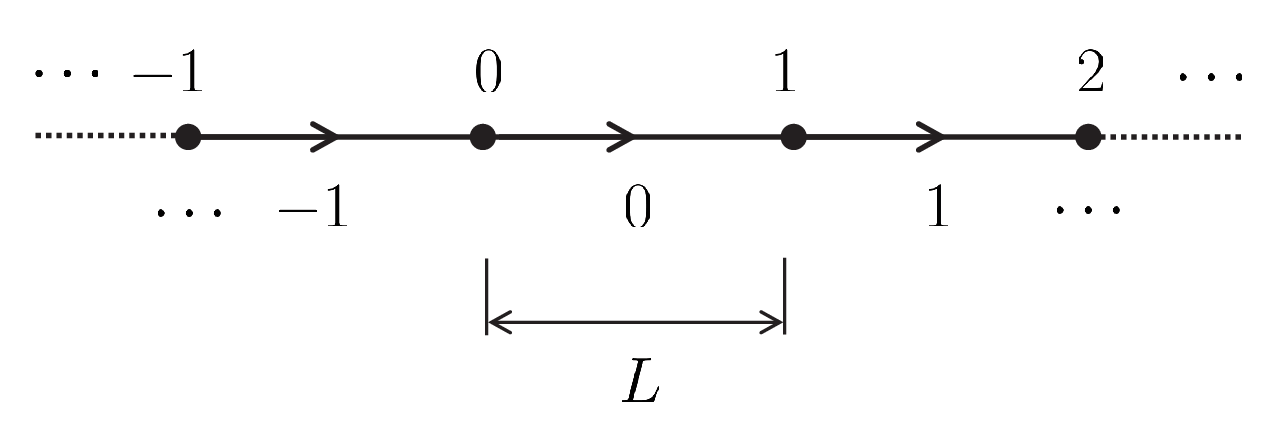}
 \end{center}
\caption{One-dimensional metamaterial. Geometry and indexing scheme for joints (top) and members (bottom). }
\label{HawmnZ}
\end{figure}

A one-dimensional metamaterial undergoing pure stretching supplies a simple example that illustrates the use of the Fourier transform to compute effective continuum properties. In this case, $n=1$, $N=1$, and the joints of the metamaterial define a simple one-dimensional Bravais lattice
\begin{equation}
  x(l) = l \, L \in \mathbb{R},
  \quad
  l \in \mathbb{Z} ,
\end{equation}
where $L$ is the length of the members, Fig.~\ref{HawmnZ}, and the discrete energy (\ref{JqkpoV}) reduces to
\begin{equation} \label{P4RLYo}
  E(u)
  =
  \sum_{m\in\mathbb{Z}} \,
    \frac{1}{2} \frac{EA}{L} (u(m+1) - u(m))^2 ,
\end{equation}
with $u(l) = v(l) \in \mathbb{R}$, $l \in \mathbb{Z}$. Alternatively, the Fourier representation of the discrete deflections is
\begin{equation} \label{Qx81lH}
    u(l)
    =
    \frac{1}{2\pi}
    \int_B
        \hat{u}(k) \, {\rm e}^{i (k L) l}
    \, dk ,
    \quad
    B = [-\frac{\pi}{L},\, \frac{\pi}{L}] ,
\end{equation}
with $\hat{u}(k) = \hat{v}(k) \in \mathbb{C}$, $k \in B$, and the dynamical matrix in
(\ref{oEf3Zx}) reduces to
\begin{equation}
  \mathbb{D}(k)
  =
  \frac{4 EA}{L^2} \, \sin^2\frac{k L}{2} ,
  \qquad
  k \in B ,
\end{equation}
which, as expected, coincides with the dynamical modulus of a harmonic monatomic chain \cite{Weiner:2002}. The scaled dynamical modulus (\ref{3ll62m}) is
\begin{equation}
  \mathbb{D}_\epsilon(k)
  =
  \frac{4 EA}{\epsilon^2 L^2} \, \sin^2\frac{\epsilon k L}{2} ,
  \qquad
  k \in B/\epsilon .
\end{equation}
Passing to the continuum limit (\ref{byu9Ue}), we obtain
\begin{equation}
  \mathbb{D}_{0}(k)
  =
  EA \, k^2 ,
  \qquad
  k \in \mathbb{R} ,
\end{equation}
which is the dynamic modulus of a linear elastic bar. The corresponding continuum energy (\ref{Lc9K2l}) is
\begin{equation} \label{DwMEu6}
    E_0(u_0)
    =
    \frac{1}{2\pi}
    \int_{\mathbb{R}}
        \frac{1}{2}
        EA \, k^2 |\hat{u}_0(k)|^2
    \, dk ,
\end{equation}
or, in real space,
\begin{equation} \label{z2tXvn}
    E_0(u_0)
    =
    \int_{\mathbb{R}}
        \frac{1}{2}
        EA \, u'^2_0(x)
    \, dx ,
\end{equation}
which is the elastic energy of a bar.

\subsection{One-dimensional metamaterial with bending} \label{1hsmET}

A one-dimensional metamaterial undergoing pure bending illustrates the scaling rule (\ref{3ll62m}). Sup\-pose that the metamaterial has the same geometry and Fourier representation as Example~\ref{TBhq2c}, Fig.~\ref{HawmnZ}, albeit with displacements $u(l) = (v(l);\, \theta(l)) \in \mathbb{R} \times \mathbb{R}$, where $v(l)$ are normal joint deflections and $\theta(l)$ are joint rotations. In this case, the discrete energy (\ref{JqkpoV}) reduces to
\begin{equation} \label{P4RLYo2}
\begin{split}
    &
    E(u)
    =
    \sum_{m\in\mathbb{Z}} \,
    \Big\{
        \frac{1}{2} \frac{EI}{L} (\theta(m+1) - \theta(m))^2
        + \\ &
        \frac{6 EI}{L}
        \Big(
            \frac{v(m+1)-v(m)}{L}
            -
            \frac{1}{2}
            \big( \theta(m) + \theta(m+1) \big)
        \Big)
    \Big\}
\end{split}
\end{equation}
or to (\ref{oEf3Zx}), in the Fourier representation, with
\begin{equation} \label{ZtWt9Q}
\begin{split}
    &
    \mathbb{D}(k)
    = \\ &
    \left(
        \begin{array}{cc}
            \dfrac{24 {EI} (1-\cos ( k L))}{L^4} &
            - \dfrac{12 i {EI} \sin ( k L)}{L^3} \\ [0.4truecm]
            \dfrac{12 i {EI} \sin ( k L)}{L^3} &
            \dfrac{4 {EI} (2 + \cos ( k L))}{L^2}
        \end{array}
    \right) .
\end{split}
\end{equation}
The scaled dynamical modulus (\ref{3ll62m}) is
\begin{equation}
\begin{split}
    &
    \mathbb{D}_\epsilon(k)
    = \\ &
    \left(
        \begin{array}{cc}
            \dfrac{24 {EI} (1-\cos ( \epsilon k L))}{\epsilon^2 L^4} &
            - \dfrac{12 i {EI} \sin ( \epsilon k L)}{\epsilon L^3} \\ [0.4truecm]
            \dfrac{12 i {EI} \sin ( \epsilon k L)}{\epsilon L^3} &
            \dfrac{4 {EI} (2 + \cos ( \epsilon k L))}{L^2}
        \end{array}
    \right) .
\end{split}
\end{equation}
Passing to the continuum limit (\ref{byu9Ue}), we obtain
\begin{equation}
\begin{split}
    &
    \mathbb{D}_0(k)
    = 
    \left(
        \begin{array}{cc}
            \dfrac{12 {EI} k^2}{L^2} &
            - \dfrac{12 i {EI} k}{L^2} \\ [0.4truecm]
            \dfrac{12 i {EI} k}{L^2} &
            \dfrac{12 {EI} }{L^2}
        \end{array}
    \right) ,
\end{split}
\end{equation}
which is the dynamic matrix of a linear elastic beam. The corresponding continuum energy (\ref{Lc9K2l}) is
\begin{equation}
    E_0(u_0)
    =
    \frac{1}{2\pi}
    \int_{\mathbb{R}}
        \frac{6 EI }{L^2}
        \, |\hat{\theta}_0(k)-i k \hat{v}_0(k)|^2
    \, dk ,
\end{equation}
or, in real space,
\begin{equation}
    E_0(u_0)
    =
    \int_{\mathbb{R}}
        \frac{6 EI }{L^2}
        \, |\theta_0(x)-v'_0(x)|^2
    \, dx ,
\end{equation}
which penalizes departures between the rotation field $\theta_0(x)$ and the deflection-compatible rotations $v'_0(x)$.

As expected, the bending strains $\theta'_0(x)$ are negligible to zeroth-order and drop out from the energy in the limit. We note that, in this example there is no axial stiffness to stabilize the continuum limit. Under these conditions, to lowest-order any pair of functions $v_0(x)$ and $\theta_0(x)$ such that $\theta_0(x)=v'_0(x)$ define a zero-energy mode and the continuum limit is unstable.

\begin{figure}[h]
 \begin{center}
 \includegraphics[width=0.45\textwidth]{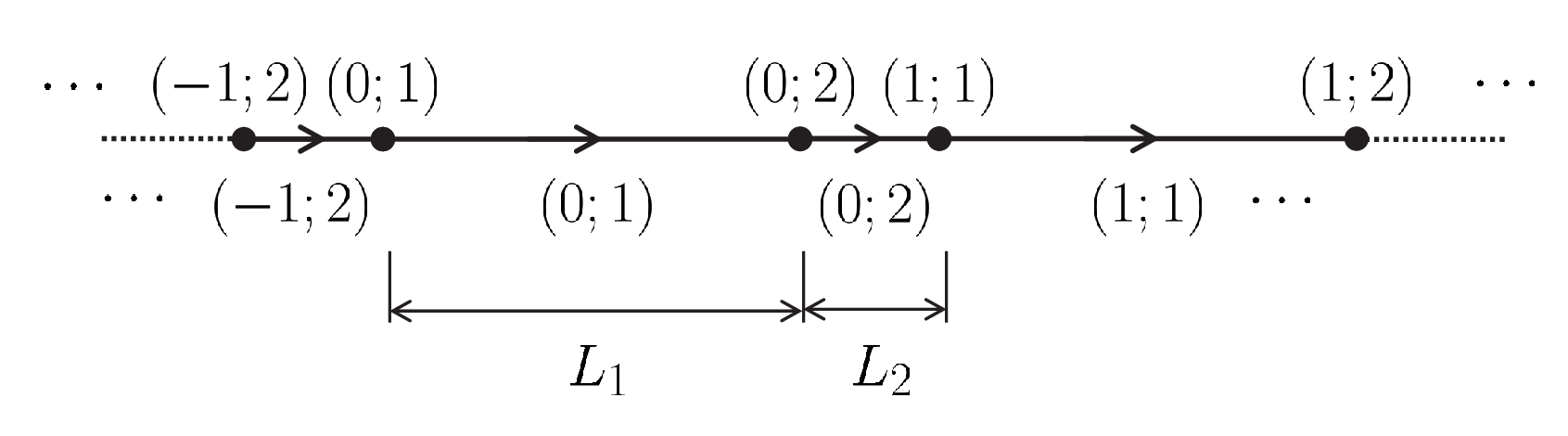}
 \end{center}
\caption{One-dimensional metamaterial with two different classes of members. Geometry and indexing scheme for joints (top) and members (bottom). }
\label{86pmH4}
\end{figure}

\subsection{Two-bar one-dimensional metamaterial}\label{ocMBEP}
A one-dimensional metamaterial combining two types of bars illustrates the energy relaxation that takes place at the unit cell level in passing to the continuum limit. In this case, $n=1$, $N=2$, and the joints of each class define the simple one-dimensional Bravais lattices
\begin{subequations}
\begin{align}
    &
    x(l;1) = l \, (L_1+L_2) ,
    &
    l \in \mathbb{Z} ,
    \\ &
    x(l;2) = l \, (L_1+L_2) + L_1 ,
    &
    l \in \mathbb{Z} ,
\end{align}
\end{subequations}
where $L_1$ and $L_2$ are the lengths of the two classes of members, cf.~Fig.~\ref{86pmH4}. Since the members can be classified into two classes of identical members, we additionally have $M=2$. The discrete energy (\ref{JqkpoV}) reduces to
\begin{equation} \label{aL7flm}
\begin{split}
    E(u)
    & =
    \sum_{m\in\mathbb{Z}}
        \frac{1}{2} \frac{EA_1}{L_1} \big(u(m;2) - u(m;1)\big)^2
    \\ & +
    \sum_{m\in\mathbb{Z}}
        \frac{1}{2} \frac{EA_2}{L_2} \big(u(m+1;1) - u(m;2)\big)^2 .
\end{split}
\end{equation}
with $u(l,\alpha) = v(l,\alpha) \in \mathbb{R}^2$, $l \in \mathbb{Z}$, $\alpha=1,2$. Alternatively, the Fourier representation of the discrete deflections is
\begin{equation}
\begin{split}
    &
    u(l,\alpha)
    =
    \frac{1}{2\pi}
    \int_B
        \hat{u}(k,\alpha) {\rm e}^{i k x(l,\alpha)}
    \, dk ,
    \\ &
    B = [-\frac{\pi}{L_1+L_2},\, \frac{\pi}{L_1+L_2}] ,
    \quad
    \alpha = 1,2,
\end{split}
\end{equation}
with $\hat{u}(k,\alpha) = \hat{v}(k,\alpha) \in \mathbb{C}^2$, $k \in B$, $\alpha=1,2$,
and the dynamical matrix in (\ref{oEf3Zx}) reduces to
\begin{equation} \label{nHrliH}
\begin{split}
    \mathbb{D}(k)
    & =
    \frac{EA_1}{L_1(L_1+L_2)}
    \left(
        \begin{array}{cc}
            1 & - {\rm e}^{- i k L_1} \\
            - {\rm e}^{  i k L_1} & 1
        \end{array}
    \right)
    \\ & +
    \frac{EA_2}{L_2(L_1+L_2)}
    \left(
        \begin{array}{cc}
            1 & -{\rm e}^{  i k L_2} \\
            -{\rm e}^{- i k L_2} & 1
        \end{array}
    \right) ,
\end{split}
\end{equation}
The scaled dynamical modulus (\ref{3ll62m}) in turn is
\begin{equation} \label{ky2nA9}
\begin{split}
    \mathbb{D}_\epsilon(k)
    & =
    \frac{EA_1}{\epsilon^2 L_1(L_1+L_2)}
    \left(
        \begin{array}{cc}
            1 & - {\rm e}^{- i \epsilon k L_1} \\
            - {\rm e}^{  i \epsilon k L_1} & 1
        \end{array}
    \right)
    \\ & +
    \frac{EA_2}{\epsilon^2 L_2(L_1+L_2)}
    \left(
        \begin{array}{cc}
            1 & -{\rm e}^{  i \epsilon k L_2} \\
            -{\rm e}^{- i \epsilon k L_2} & 1
        \end{array}
    \right) ,
\end{split}
\end{equation}
Passing to the continuum limit (\ref{U2E0YW}), with
\begin{equation}
    L = \Big( \frac{1}{2},\, \frac{1}{2} \Big) ,
\end{equation}
gives
\begin{equation}
    \mathbb{D}_0(k)
    =
    \Big(
        \frac{L_1}{L_1+L_2}
        \frac{1}{EA_1}
        +
        \frac{L_2}{L_1+L_2}
        \frac{1}{EA_2}
    \Big)^{-1}
    \, k^2 ,
\end{equation}
where, in the prefactor of $k^2$, we recognize the effective modulus $EA$ of two bars connected in series, as expected. With this identification, the continuum energy reduces to the forms (\ref{DwMEu6}) and (\ref{z2tXvn}).

It bears emphasis that the attainment of the correct limit requires the cell-wise relaxation step encoded in the limit (\ref{U2E0YW}). Indeed, we note from (\ref{ky2nA9}) that $\mathbb{D}_\epsilon(k)$ diverges as $\epsilon\to 0$ and, therefore, a na\"ieve pointwise limit of the energy fails. We also note that, contrariwise, $\lim_{\epsilon\to 0} \mathbb{D}_\epsilon^{-1}(k)$ is finite, but the limiting matrix is singular for fixed $k$ and cannot be inverted. Only when $\lim_{\epsilon\to 0} \mathbb{D}_\epsilon^{-1}(k)$ is projected by $L$, as in (\ref{U2E0YW}), is the resulting matrix $\mathbb{D}_0^{-1}(k)$ invertible and a proper limit $\mathbb{D}_0(k)$ is defined.

\subsection{Two-dimensional honeycomb metamaterial}\label{xFzCtx}

\begin{figure*}[ht!]
\begin{center}
    \subfigure[]{\includegraphics[width=0.4\textwidth]{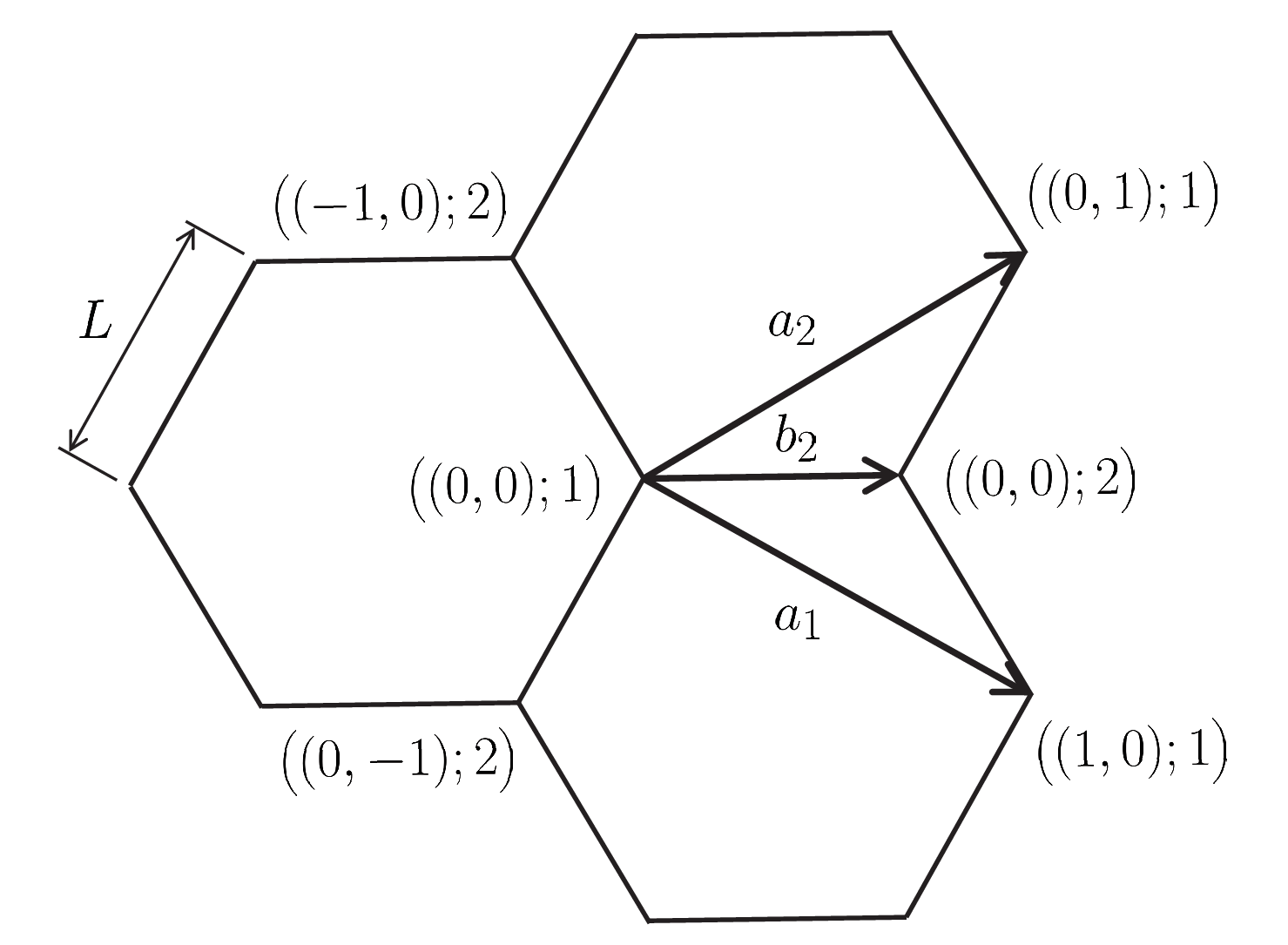}}
    \subfigure[]{\includegraphics[width=0.32\textwidth]{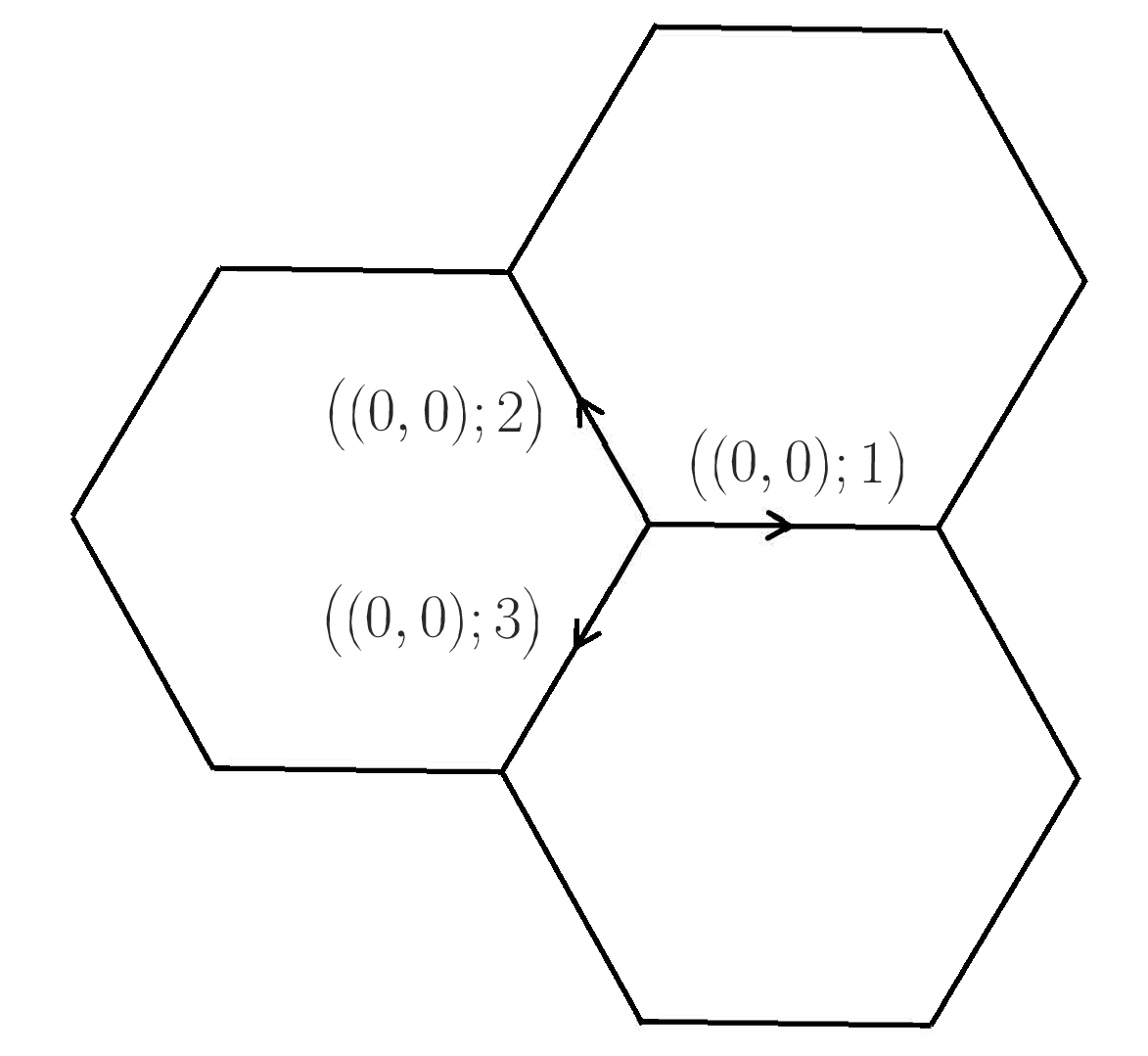}}
    \caption{Honeycomb metamaterial. (a) Joint numbering scheme using two simple Bravais sublattices. (b) Member numbering scheme using three simple Bravais sublattices.} \label{FIDpWS}
\end{center}
\end{figure*}

The honeycomb metamaterial of size $L$ is two dimensional, $n=2$, and contains two types of joints, $N=2$, and three types of oriented bars, $M=3$, Fig.~\ref{FIDpWS}. A possible choice of Bravais lattice and shifts is, Fig.~\ref{FIDpWS}a,
\begin{subequations}
\begin{align}
    &
    a_1 = L \, ( \frac{3}{2},\, - \frac{\sqrt{3}}{2} ) ,
    &
    b_1 = L \, (0,\, 0) ,
    \\ &
    a_2 = L \, ( \frac{3}{2},\, + \frac{\sqrt{3}}{2} ) ,
    &
    b_2 = L \, (1,\, 0) .
\end{align}
\end{subequations}
The volume of the corresponding unit cell is \cite{Ariza:2010}
\begin{equation}
    V = \frac{ 3 \sqrt{3}}{2} \, L^2 .
\end{equation}
Fig.~\ref{FIDpWS}b shows a choice of indexing for the three types of bars and a choice of bar orientation. The spanning vectors for the three classes of bars are
\begin{equation}
\begin{split}
    &
    dx_1 = \{1,\, 0\} \, L,
    \\ &
    dx_2 = \{-\frac{1}{2},\, \frac{\sqrt{3}}{2}\} \, L ,
    \\ &
    dx_3 = \{-\frac{1}{2},\, -\frac{\sqrt{3}}{2} \} \, L .
\end{split}
\end{equation}
A direct evaluation of (\ref{59joLe}), (\ref{XLBfZ6}) and (\ref{OVv27j}), followed by a computation of the limit (\ref{U2E0YW}), gives the limiting continuum energy
\begin{equation}
\begin{split}
    &
    E_0(u_0)
    =
    \int_\Omega
        \frac{1}{2}
        \mathbb{C} \,
        \varepsilon(x)
        \cdot
        \varepsilon(x)
    \, dx
    + \\ &
    \int_\Omega
        \frac{1}{2}
        \frac{8\sqrt{3} EI}{L}
        \Big(
            \theta(x)
            -
            \frac{1}{2}
            \big(
                \partial_1 v_2(x)
                -
                \partial_2 v_1(x)
            \big)
        \Big)^2
    \, dx
\end{split}
\end{equation}
in terms of elastic and micropolar components, respectively. In the elastic part of the energy,
\begin{equation}
    \varepsilon
    =
    \left\{
        \begin{array}{c}
            \varepsilon_{11} \\
            \varepsilon_{22} \\
            2 \varepsilon_{12}
        \end{array}
    \right\}
    =
    \left\{
        \begin{array}{c}
            \partial_1 v_1 \\
            \partial_2 v_2 \\
            \partial_1 v_2 + \partial_2 v_1
        \end{array}
    \right\}
\end{equation}
are continuum strains and
\begin{equation}
    \mathbb{C}
    =
    \left(
        \begin{array}{ccc}
            C_{11} & C_{12} & C_{13} \\
            C_{12} & C_{22} & C_{23} \\
            C_{13} & C_{23} & C_{33}
        \end{array}
    \right)
\end{equation}
are effective elastic moduli, explicitly,
\begin{equation} \label{YRC3k2}
\begin{split}
    &
    C_{11}
    =
    C_{22}
    =
    \frac
    {
        {EA} \, \left(\, {EA} \, L^2 + 36 \, {EI} \, \right)
    }
    {
        2 \sqrt{3} \, \left(\, {EA} \, L^3 + 12 \, {EI} \, L\right)
    } ,
    \\ &
    C_{33}
    =
    \frac
    {
        4 \sqrt{3} \, {EA} \, {EI}
    }
    {
        {EA} \, L^3 + 12 \, {EI} \, L
    } ,
    \\ &
    C_{12}
    =
    \frac
    {
        {EA} \, \left({EA} \, L^2 - 12 \, {EI} \, \right)
    }
    {
        2 \sqrt{3} \,
        \left(\, {EA} \, L^3 + 12 \, {EI} \, L\right)
    } ,
    \\ &
    C_{13} = C_{23} = 0 .
\end{split}
\end{equation}
We note that the continuum energy is isotropic in the plane, as expected from the hexagonal symmetry of the lattice.  We also note the rational dependence of the elastic moduli (\ref{YRC3k2}) on bar properties, which is the result of local relaxation at the cell level between the two types of nodes in the honeycomb structure. Finally, we find that the elastic moduli are singular if $EI=0$, i.~e., for the articulated honeycomb truss. This is indeed expected, since the honeycomb truss is a mechanism with zero-energy modes of deformation. The honeycomb structure thus illustrates the critical role of bending in stabilizing bending-dominated metamaterials.

\subsection{Three-dimensional octet metamaterial}\label{obmfAG}

\begin{figure*}[h!]
\begin{center}
    \subfigure[]{\includegraphics[width=0.7\textwidth]{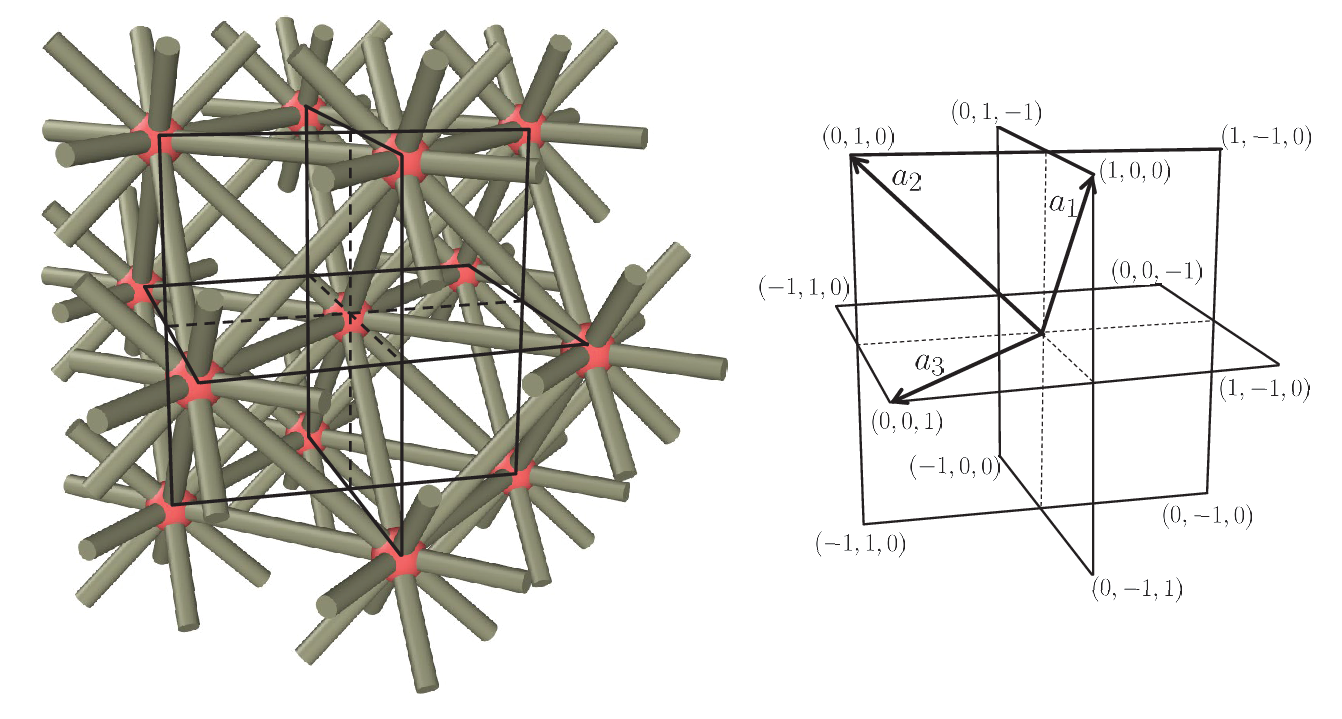}}
    \subfigure[]{\includegraphics[width=0.6\textwidth]{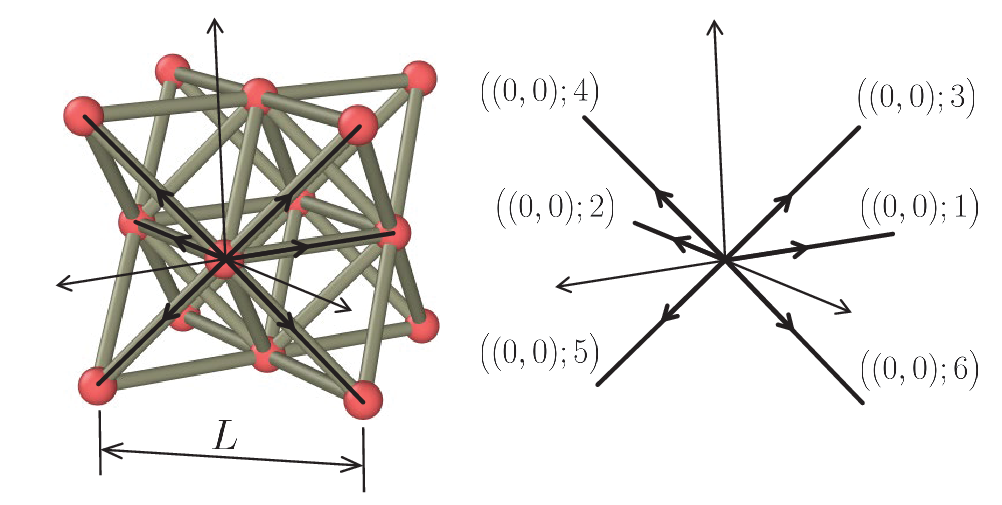}}
    \caption{Octet-truss metamaterial. (a) Joint numbering scheme using a simple Bravais lattice. (b) Bar numbering scheme.} \label{S2rFkk}
\end{center}
\end{figure*}

The octet-truss metamaterial of size $L$ is three dimensional, $n=3$, and contains one type of joints, $N=1$, and six types of oriented bars of length $L/\sqrt{2}$, $M=6$, Fig.~\ref{S2rFkk}. A possible choice of Bravais lattice is, Fig.~\ref{S2rFkk}a,
\begin{subequations}
\begin{align}
    &
    a_1 = \frac{L}{\sqrt{2}} \, \{ 0,\, 1,\, 1
    \} ,
    \\ &
    a_2 = \frac{L}{\sqrt{2}} \, \{ 1,\, 0,\, 1 \},
    \\ &
    a_3 = \frac{L}{\sqrt{2}} \, \{ 1,\, 1,\, 0 \} .
\end{align}
\end{subequations}
In addition, Fig.~\ref{S2rFkk}b shows an indexing for the three types of bars and a choice of bar orientation. The volume of the
corresponding unit cell is \cite{Ariza:2005}
\begin{equation}
    V = \frac{L^3}{\sqrt{2}} .
\end{equation}
The bar vectors for the six classes of bars are
\begin{equation}
\begin{split}
    &
    dx_1 = - a_3,
    \quad
    dx_2 = - a_1 + a_2 ,
    \quad
    dx_3 = a_1 ,
    \\ &
    dx_4 = a_2 ,
    \quad
    dx_5 = - a_1 + a_3 ,
    \quad
    dx_6 = - a_2 + a_3 .
\end{split}
\end{equation}
For simplicity, we assume isotropic moments of inertia, $I_2 = I_3 := I$, e.~g., corresponding to a thin-walled circular cross section. Under this assumption, the choice of directors $d_{2\beta}$ and $d_{3\beta}$ is immaterial, provided that they define an orthonormal triad together with $d_{1\beta}$. For definiteness, in calculations we choose
\begin{equation}
    d_{1\beta} = \frac{dx_\beta}{|dx_\beta|} ,
    \;\;
    d_{2\beta} = \frac{d_{1\beta} \times e_3}{|dx_\beta \times e_3|} ,
    \;\;
    d_{3\beta} = d_{1\beta} \times d_{2\beta} ,
\end{equation}
where $e_3 = (0,0,1)$ in the Cartesian reference frame.

A direct evaluation of (\ref{59joLe}), (\ref{XLBfZ6}) and (\ref{OVv27j}), followed by a computation of the limit (\ref{U2E0YW}), gives the limiting continuum energy
\begin{equation}
\begin{split}
    &
    E_0(u_0)
    =
    \int_\Omega
        \frac{1}{2}
        \mathbb{C} \,
        \varepsilon(x)
        \cdot
        \varepsilon(x)
    \, dx
    +
    \int_\Omega
        \frac{1}{2}
        \frac{48 \sqrt{2} {EI}}{L^4}
        \Big(
        \\ & \quad\quad\quad
            \big(
                \theta_1(x)
                -
                \frac{1}{2}
                (\partial_2 v_3(x) - \partial_3 v_2(x) )
            \big)^2
            + \\ & \quad\quad\quad
            \big(
                \theta_2(x)
                -
                \frac{1}{2}
                (\partial_3 v_1(x) - \partial_1 v_3(x) )
            \big)^2
            + \\ & \quad\quad\quad
            \big(
                \theta_3(x)
                -
                \frac{1}{2}
                (\partial_1 v_2(x) - \partial_2 v_1(x) )
            \big)^2
        \Big)
    \, dx
\end{split}
\end{equation}
in terms of elastic and micropolar components, respectively. In the elastic part of the energy,
\begin{equation}
    \varepsilon
    =
    \left\{
        \begin{array}{c}
            \varepsilon_{11} \\
            \varepsilon_{22} \\
            \varepsilon_{33} \\
            2 \varepsilon_{12} \\
            2 \varepsilon_{13} \\
            2 \varepsilon_{23}
        \end{array}
    \right\}
    =
    \left\{
        \begin{array}{c}
            \partial_1 v_1 \\
            \partial_2 v_2 \\
            \partial_3 v_3 \\
            \partial_1 v_2 + \partial_2 v_1 \\
            \partial_1 v_3 + \partial_3 v_1 \\
            \partial_2 v_3 + \partial_3 v_2
        \end{array}
    \right\}
\end{equation}
are continuum strains and
\begin{equation}
    \mathbb{C}
    =
    \left(
        \begin{array}{cccccc}
            C_{11} & C_{12} & C_{12} & 0 & 0 & 0 \\
            C_{12} & C_{11} & C_{12} & 0 & 0 & 0 \\
            C_{12} & C_{12} & C_{11} & 0 & 0 & 0 \\
            0 & 0 & 0 & C_{44} & 0 & 0 \\
            0 & 0 & 0 & 0 & C_{44} & 0 \\
            0 & 0 & 0 & 0 & 0 & C_{44}
        \end{array}
    \right)
\end{equation}
are effective elastic moduli, explicitly,
\begin{equation}
\begin{split}
    &
    C_{11}
    =
    \frac{4 {EA} L^2+24 {EI}}{\sqrt{2} L^4} ,
    \\ &
    C_{12}
    =
    \frac{\sqrt{2} \left({EA} L^2-6 {EI}\right)}{L^4} ,
    \\ &
    C_{44}
    =
    \frac{2 {EA} L^2+12 {EI}}{\sqrt{2} L^4} .
\end{split}
\end{equation}
We verify that $C_{11} \neq C_{12} + 2 C_{44}$. Hence, the continuum strain-energy density has cubic symmetry, while the micropolar energy density is isotropic, as expected from the cubic symmetry of the lattice. It is also interesting to note that the effective moduli are independent of the torsional stiffness $GI_1$ of the bars.

\section{Summary and discussion}

By an appeal to concepts of calculus of variations, including notions of variational convergence ($\Gamma$-convergence in mathematical parlance), we have elucidated the continuum limit of linear metastructures under the assumption of a strict separation of micromechanical (small cell size) and macromechanical (large structure) length scales. The continuum limit of metastructures is of the discrete-to-continuum type \cite{Alicandro:2004,BraidesGelli2006}, but it is non-standard due to presence of bending and rotational degrees of freedom. By an interpolation argument, the continuum limit of metastructures can also be formulated as a homogenization limit \cite{Cioranescu:1999}. Based on this connection, we argue that the limiting continuum energy is of the integral type, with a well-defined energy density that depends on local gradients and attendant natural boundary conditions. The elucidation of the limiting energy density is greatly facilitated by the assumption of linear elasticity and the corresponding quadratic form of the discrete energy. We have presented a procedure based on the discrete Fourier transform that yields explicitly the effective moduli in the continuum limit, and illustrated the procedure by means of elementary examples and two common metamaterials: a two-dimensional honeycomb lattice and a three-dimensional octet-truss.

The analysis presented in the foregoing sheds light on---and raises---a number of questions, which we briefly discuss next in closing.

\underline{\bf The variational continuum limit is unique.} The \newline variational limit of an infinitely fine metastructure of given domain $\Omega$ is characterized by the unique continuum energy that matches exactly the limiting minimum potential energies of a sequence of increasingly fine metastructures {\sl for all admissible forcings} \cite{dalmaso:1993}. By extension, it is also the unique continuum energy that matches exactly the limiting values of any macrostructural quantity of interest that depends continuously on the total energy, such as reactions or load-displacement relations. In this sense, the variational continuum limit fulfills its intended aim of being {\sl indistinguishable} to a macroscopic observer from a sufficiently fine metastructure.

The variational limit also has the unique property of being {\sl lossless}, in the sense that any displacement solution of the limiting continuum problem is the limit, in a precise sense, of the discrete displacement solutions of the corresponding sequence of increasingly fine metastructures. In particular, the approximating sequence of discrete displacements can be {\sl reconstructed}, is desired, from the continuum solution at {\sl no loss of information}.

By virtue of these properties, the variational continuum limit is optimal and superior to any other {\sl ad hoc} continuum model thereof. In particular, {\sl ad hoc} continuum models of fine metastructures different from the variational continuum limit are sure to result in {\sl gaps} between the continuum energy and the limiting energies of the metastructures, and, by extension, also in gaps in any other macroscopic quantity of interest, at least for some forcings.

\underline{\bf The continuum limit is micropolar and local.} By a rigorous computation of the variational continuum limit of infinitely fine metastructures, we have shown that the limiting continuum energy is {\sl micropolar}, in the sense of \cite{Eringen:1966}. However, to lowest order the continuum energy is independent of the {\sl curvature} of the deformation, or {\sl bending strain}, and therefore {\sl local}. In particular, the lowest-order continuum limit exhibits linear-elastic volume scaling and fails to capture {\sl size effects} such as observed in fracture experiments at the mesoscale \cite{Shaikeea2022}. The scope of the lowest-order continuum limit is therefore limited to applications in which size effects can be neglected, such as the analysis of lightweight structures under normal operating conditions, cf., e.~g., \cite{Cheung:2013, Gregg:2024} for representative applications to aerospace structures.

\underline{\bf Higher-order expansions and size effects.} Ev\-i\-dent\-ly, the local character of the continuum limit to lowest order owes to the assumption of strict scale separation. In situations where non-local or mesoscale effects are important, the lowest-order limit needs to be augmented with fine-scale information. The full micropolar model of the honeycomb metastructure proposed by \cite{CHEN1998} is an early example in that regard. As mentioned in the introduction, computational schemes have been proposed that coarse-grain the description of a metastructure while retaining full resolution where required \cite{kochmann:2019, PHLIPOT2019}.

Alternatively, higher-order generalized continuum models can be systematically formulated by recourse to the method of variational $\Gamma$-equivalence \cite{Braides:2008}. We recall that a sequence $(G_\epsilon)$ of energy functionals is {\sl variationally equivalent} to another $(F_\epsilon)$ to order $\alpha > 0$ as $\epsilon \to 0$ if
\begin{equation}
  m(F_\epsilon; f_0) \sim m(G_\epsilon; f_0) + o(\epsilon^\alpha) ,
\end{equation}
{\sl for all admissible forcings} $f_0$, where $m(F_\epsilon; f_0)$, respectively $m(G_\epsilon; f_0)$, is the minimum potential attained by $F_\epsilon$, respectively $G_\epsilon$, under forcing $f_0$. The fundamental property of variational equivalence (see \cite[Theorems 4.3 and 4.4]{Braides:2008} is that, provided that $(F_\epsilon)$ and $(G_\epsilon)$ are stable in the sense of equicoercivity, their low-energy configurations converge to the same limit. In addition, any macroscopic quantity of interest computed from $(G_\epsilon)$ is indistinguishable from the same quantity computed from $(F_\epsilon)$ to order $o(\epsilon^\alpha)$.

The variational $\Gamma$-equivalence framework thus provides a rigorous means of defining sequences of mesoscopic continuum models $(G_\epsilon)$ that have the same continuum limit $F_0$ as the sequence of discrete metastructural energies $(F_\epsilon)$ and, in addition, exhibit the same asymptotic behavior to order $o(\epsilon^\alpha)$. For instance, for $N=1$ it is readily verified that a formal Taylor series expansion of $\mathbb{D}^{-1}_\epsilon(k)$ to order $\alpha$ in $\epsilon$, namely,
\begin{equation}\label{L66qmN}
  \mathbb{D}_{\epsilon,\alpha}(k)
  :=
  \Big(
    \sum_{m=0}^\alpha
    \frac{1}{m!}
    \Big(
      \frac{d^m}{d\epsilon^m}
      \mathbb{D}_\epsilon^{-1}(k)
    \Big)_{\epsilon=0}
    \, \epsilon^m
  \Big)^{-1} ,
\end{equation}
generates a sequence of continuum energies $(F_{\epsilon,\alpha})$ that is variationally equivalent to the discrete metastructural energies $(F_\epsilon)$ to order $\alpha$. It should be carefully noted that variationally equivalent energies of the same order are not unique, and particular choices thereof are, to a certain extent, a matter of convenience or expediency. However, variational equivalence does provide a rigorous criterion to weed out non-equivalent and, therefore, non-performing models.

In the present context, the higher-order extensions of the continuum limit, such as the full micropolar model proposed by \cite{CHEN1998} for honeycomb metastructures, are non-local and therefore exhibit a size effect. By assuming the bars to have a finite strength, they also endow the limiting continuum with a well-defined fracture toughness \cite{CHEN1998}. It remains to verify whether these enhancements of the theory suffice to capture observed {\sl size effects}, such as the specimen-size dependence of fracture toughness \cite{Shaikeea2022}.

\section*{Acknowledgements}

S.~Conti and M.~Ortiz gratefully acknowledge the support of the Deutsche Forschungsgemeinschaft (DFG, German Research Foundation) {\sl via} project 211504053 - SFB 1060; project 441211072 - SPP 2256; and project 390685813 -  GZ 2047/1 - HCM. P.~Ariza gratefully acknowledges financial support from Ministerio de Ciencia e Innovación under grant number PID2021-124869NB-I00.

\begin{appendix}

\section{The discrete Fourier transform} \label{TwWqmj}

Let $(a_i)_{i=1}^n$ be a basis of $\mathbb{R}^n$ and $\mathcal{L} = \{ x(l) := \sum_{i=1}^n l^i a_i \, : \, l \in \mathbb{Z}^n \}$ the corresponding Bravais lattice. Let $f : \mathcal{L} \to \mathbb{R}$ be a real-valued lattice function. The discrete Fourier transform of $f$ is a complex function $\hat{f}(k)$ supported on the Brillouin zone $B$ in dual space given by
\begin{equation} \label{IYOTsT}
  \hat{f}(k)
  =
  V \sum_{l \in \mathbb{Z}^n}
  f(l) {\rm e}^{-i k \cdot x(l) } .
\end{equation}
where $V$ is the volume of the unit cell of the lattice. The inverse mapping is given by
\begin{equation} \label{aX0TEN}
  f(l)
  =
  \frac{1}{(2\pi)^n} \int_B
  \hat{f}(k) {\rm e}^{i k \cdot x(l) } dk .
\end{equation}
The convolution of two lattice functions $f(l)$, $g(l)$ is
\begin{equation} \label{YNU3CC}
  (f*g)(l) = V \sum_{l' \in \mathbb{Z}^n} f(l-l') g(l') ,
\end{equation}
whereupon the convolution theorem states that
\begin{equation} \label{m4IeYG}
  \widehat{f*g} = \hat{f} \hat{g} .
\end{equation}
In addition, the Parseval identity states that
\begin{equation} \label{pB0rWJ}
  V \sum_{l \in \mathbb{Z}^n} f(l) g^*(l)
  =
  \frac{1}{(2\pi)^n} \int_{B}
  \hat{f}(k) \hat{g}^*(k) dk
\end{equation}
which establishes an isometric isomorphism between $l^2$ and $L^2(B)$ (cf.~\cite{Ariza:2005} and references therein for further details).

\section{Structure of continuum energy} \label{KA4bct}

Let $\lambda > 0$. From (\ref{3ll62m}), we have
\begin{equation} \label{RWPF24}
    \mathbb{D}_\epsilon(\lambda k)
    \, \zeta \cdot \zeta^*
    =
    \epsilon^{-2} \mathbb{D}(\epsilon \lambda k)
    (\xi; \, \epsilon \eta)
    \cdot
    (\xi^*; \, \epsilon \eta^*) .
\end{equation}
Rescale $\epsilon$ as
\begin{equation}
    \delta = \lambda \epsilon,
    \quad
    \epsilon = \lambda^{-1} \delta .
\end{equation}
Then, from (\ref{RWPF24}),
\begin{equation}
\begin{split}
    &
    \mathbb{D}_\epsilon(\lambda k)
    \, \zeta \cdot \zeta^*
    = \\ &
    \lambda^2 \delta^{-2} \mathbb{D}(\delta k)
    \,
    (\xi;\, \lambda^{-1} \delta \eta)
    \cdot
    (\xi^*;\, \lambda^{-1} \delta \eta^*)
    = \\ &
    \lambda^2 \mathbb{D}_\delta(k)
    \,
    (\xi;\, \lambda^{-1} \eta)
    \cdot
    (\xi^*;\, \lambda^{-1} \eta^*) ,
\end{split}
\end{equation}
where we have used (\ref{3ll62m}) again. By duality,
\begin{equation} \label{Vr5RU0}
\begin{split}
    &
    \frac{1}{2}
    \mathbb{D}_\epsilon^{-1}(\lambda k)
    \, {\gamma} \cdot \gamma^*
    = \\ &
    \min_{\zeta}
    \{
        \gamma \cdot \zeta^*
        -
        \frac{1}{2}
        \lambda^2 \mathbb{D}_\delta(k)
        \,
        (\xi;\, \lambda^{-1} \eta)
        \cdot
        (\xi^*;\, \lambda^{-1} \eta^*)
    \} .
\end{split}
\end{equation}
Writing
\begin{equation}
    \gamma := (\alpha; \, \beta)
    \in [\mathbb{C}^n]^N \times [\mathbb{C}^{n(n-1)/2}]^N ,
\end{equation}
we have the decomposition
\begin{equation}
    \gamma \cdot \zeta^*
    =
    \alpha \cdot \xi^* + \beta \cdot \eta^*
    =
    (\alpha;\, \lambda \beta)
    \cdot
    (\xi^*;\, \lambda^{-1} \eta^*) ,
\end{equation}
in terms of deflections and rotations. Inserting into (\ref{Vr5RU0}) we obtain, explicitly,
\begin{equation}
    \mathbb{D}_\epsilon^{-1}(\lambda k)
    \, {\gamma} \cdot \gamma^*
    =
    \lambda^{-2} \mathbb{D}_\delta^{-1}(k) \,
    (\alpha;\, \lambda \beta)
    \cdot
    (\alpha^*;\, \lambda \beta^*) .
\end{equation}
Passing to the limit using (\ref{U2E0YW}) further gives
\begin{equation}
\begin{split}
    &
    \mathbb{D}_0^{-1}(\lambda k)
    \, {\gamma}_0 \cdot {\gamma}_0^*
    = \\ &
    \lim_{\delta\to 0}
    \lambda^{-2} \mathbb{D}_\delta^{-1}(k) \,
    L(\alpha_0;\, \lambda \beta_0)
    \cdot
    L(\alpha_0^*;\, \lambda \beta_0^*)
    = \\ &
    \lambda^{-2} \mathbb{D}_0^{-1}(k) \,
    (\alpha_0;\, \lambda \beta_0)
    \cdot
    (\alpha_0^*;\, \lambda \beta_0^*) ,
\end{split}
\end{equation}
with
\begin{equation} \label{cSASuq}
    \gamma_0 := (\alpha_0; \, \beta_0)
    \in \mathbb{C}^n \times \mathbb{C}^{n(n-1)/2} .
\end{equation}
Finally, by duality we have
\begin{equation} \label{5elpc0}
\begin{split}
    &
    \frac{1}{2}
    \mathbb{D}_0(\lambda k)
    \, {\zeta}_0 \cdot {\zeta}_0^*
    =
    \min_{\gamma_0}
    \{
        \gamma_0 \cdot {\zeta}_0^*
        - \\ &
        \frac{1}{2}
        \lambda^{-2} \mathbb{D}_0^{-1}(k) \,
        (\alpha_0;\, \lambda \beta_0)
        \cdot
        (\alpha_0^*;\, \lambda \beta_0^*)
    \} .
\end{split}
\end{equation}
Decomposing $\zeta_0$ as in (\ref{8t4tSW}), to match (\ref{cSASuq}), the duality pairing has the decomposition
\begin{equation}
    \gamma_0 \cdot \zeta_0^*
    =
    \alpha_0 \cdot \xi_0^* + \beta_0 \cdot \eta_0^*
    =
    (\alpha_0;\, \lambda \beta_0)
    \cdot
    (\xi_0^*;\, \lambda^{-1} \eta_0^*) ,
\end{equation}
and (\ref{5elpc0}) evaluates to
\begin{equation}
\begin{split}
    &
    \mathbb{D}_0(\lambda k)
    \, {\zeta}_0 \cdot {\zeta}_0^*
    = \\ &
    \lambda^2
    \mathbb{D}_0(k) \,
    (\xi_0;\, \lambda^{-1} \eta_0)
    \cdot
    (\xi_0^*;\, \lambda^{-1} \eta_0^*) ,
\end{split}
\end{equation}
or equivalently (\ref{qwzrN3}), as surmised.

\end{appendix}

\bibliography{biblio}
\bibliographystyle{unsrt}

\end{document}